\documentclass[11pt]{amsart}
\usepackage{amsmath,amssymb,amscd}
\usepackage{mathrsfs}
\usepackage{graphicx}
\usepackage[all]{xy}

\title{Minitwistor spaces, Severi varieties, and Einstein-Weyl structure}

\author{Nobuhiro Honda, Fuminori Nakata}
\thanks
{Both authors are partially supported by
Research Fellowships of the 
Japan Society for the Promotion
of Science for Young Scientists.
}
\date{}

\newcommand{\ol}{\overline}

\newcommand{\lra}{\longrightarrow}
\newcommand{\lras}{\,\longrightarrow\,}

\newcommand{\set}{\,|\,}
\newcommand{\proofend}{\hfill$\square$}

\newtheorem{prop}{Proposition}[section]
\newtheorem{lemma}[prop]{Lemma}

\newtheorem{thm}[prop]{Theorem}
\newtheorem{rmk}[prop]{Remark}
\newtheorem{cor}[prop]{Corollary}

\newtheorem{definition}[prop]{Definition}

\begin{document}

\begin{abstract}
In this paper we show that the space  of {\em nodal}\, rational curves, which is so called a Severi variety (of rational curves),  on any non-singular projective  surface is always  equipped with a natural Einstein-Weyl structure, if the space is 3-dimensional. This is a generalization of the Einstein-Weyl structure on the space  of {\em smooth} rational curves on a complex surface, given by N.\,Hitchin.
As geometric objects naturally associated to Einstein-Weyl structure,
we investigate null surfaces and  geodesics on the Severi varieties.
Also we see that if the projective surface has an appropriate real structure, then the real locus of the Severi variety becomes a positive definite Einstein-Weyl manifold. 
Moreover we  construct various explicit examples of rational surfaces having 3-dimensional Severi varieties of rational curves.
\end{abstract}

\maketitle
\section{Introduction}
In the paper \cite{H82-1}, N.\,J.\,Hitchin established a kind of twistor correspondence which provides a bijection between  
3-dimensional Einstein-Weyl  manifolds and  non-singular complex surfaces which have  non-singular rational curves with normal bundle $\mathscr O(2)$.
The latter complex surfaces are called the {\em minitwistor spaces}, and the rational curves in the spaces  are called the {\em minitwistor lines}.
In this Hitchin correspondence,  Einstein-Weyl manifolds appear as   parameter spaces of the minitwistor lines.
The parameter space has a natural complex conformal structure and holomorphic geodesics, by which the space becomes an Einstein-Weyl 3-fold.
Conversely, when an Einstein-Weyl 3-fold is given, the minitwistor space is obtained as the leaf space for a certain foliation on a  conic bundle naturally constructed over the Einstein-Weyl 3-fold. 
Here in order to obtain the minitwistor space, we need to assume in general that an Einstein-Weyl 3-fold is sufficiently small.
In this sense, the Hitchin correspondence is local in nature.
It is also remarkable that there are essentially only two compact minitwistor spaces,
while there are many non-compact minitwistor spaces as constructed by H.\,Pedersen and K.\,P.\,Tod (\cite{Pe86,PT93}). 
The two compact minitwistor spaces correspond to the two standard Einstein-Weyl 3-folds (i.e.\,the Euclidean space and the hyperbolic space). 

In this paper, we show that  if we allow the minitwistor lines to be {\em nodal}\, rational curves, then their parameter space still carries a  natural Einstein-Weyl structure, if the parameter space is 3-dimensional.
When the complex surface is projective algebraic, 
the parameter space of nodal curves (of any genus) in a linear system is called a {\em Severi variety} in algebraic geometry, which is known to have a natural structure of a (non-complete) algebraic variety.
Thus our result can be precisely stated that\! {\em any Severi variety of rational curves has a natural Einstein-Weyl structure, if it is 3-dimensional}.
If $C$ denotes any of the nodal rational curves, 
the last condition is equivalent to the condition that the self-intersection number of $C$ and the number of the nodes of $C$ are $2m$ and $(m-1)$ respectively. 
So when $m=1$ our result is reduced to the Hitchin's original Einstein-Weyl structure.

In Section 2.1 we recall the above Hitchin correspondence, with some emphasizes on null plane bundles over Einstein-Weyl 3-folds, a plane distribution on the bundle, and its integrability.
In Section 2.2 we first recall fundamental results on Severi varieties in general (Proposition \ref{prop-severi}), and show that if  the nodal curves (parametrized by the Severi variety) are rational, then the Severi variety is non-singular and its dimension is expressed by the self-intersection number of the rational curves and the number of the nodes (Proposition \ref{prop-severi2}).
In Section 2.3 by using Hitchin's result we show that the Severi variety of rational curves has a natural Einstein-Weyl structure, if the variety is 3-dimensional (Theorem \ref{thm:definingEWstr}).
In Section 2.4 motivated by this result we define {\em minitwistor spaces} to be a pair of non-singular  projective surface and a linear system on it which has a 3-dimensional Severi variety of rational curves (Definition \ref{def:mt}). 
Our definition involves a positive integer $m$ which is one greater than the number of the nodes, and we call this integer the {\em index} of the minitwistor space. Then the original minitwistor spaces by Hitchin are exactly the minitwistor spaces of index one.
As any blow-up of the minitwistor spaces (in the above sense) becomes again minitwistor space, we also introduce the notion of minimality of the minitwistor spaces (Definition \ref{def:minimal}). 
Then any minimal minitwistor spaces of index one are isomorphic to either $\mathbb{CP}^1\times\mathbb{CP}^1$ or the Hirzebruch surface $\mathbb P(\mathscr O(2)\oplus\mathscr O)$ (Proposition \ref{prop:class1}).

In Section 3 we investigate certain subvarieties of the 3-dimensional Severi variety (of rational curves) naturally arising from the Einstein-Weyl structure.
Namely we investigate null surfaces and (null and non-null) geodesics in the Einstein-Weyl 3-fold.
Just as in the Hitchin's case, null surfaces are formed by minitwistor lines going through a point on the minitwistor space, and geodesics are formed by those going through two points on the minitwistor space.
(In particular, they are automatically algebraic subvarieties.)
But significant difference in our case is that both of them are  non-normal subvarieties.
The singular locus of these subvarieties are formed by minitwistor lines which have nodes at the prescribed point(s).

In Section 4 we show that if the minitwistor space (in our sense) has a real structure and  a real  minitwistor line which has the nodes as its all real points, then the real locus of the (3-dimensional) Severi variety has a natural structure of a {\em real, positive definite} Einstein-Weyl 3-manifold (Theorem \ref{thm:real}).  
They are obtained 
as  real slices of the complex Einstein-Weyl structure obtained in Section 
\ref{Section:Hitchin&Severi}.

In Section 5 we provide various examples of the minitwistor spaces.
In Section 5.1 for any $m\ge 2$ we construct minimal minitwistor spaces of index $m$.
They are obtained from the product surface $\mathbb{CP}^1\times\mathbb{CP}^1$ by blowing-up $2m$ points.
 In this example, the configuration of the  $2m$ points can be taken generically, so that they constitute a $4m$-dimensional family, while  the number of effective parameter is $4m-6$.
Also we see that if we specialize the configuration of the  $2m$ points  in a certain way, then we obtain minitwistor spaces with $\mathbb C^*$-action, or even {\em toric} minitwistor spaces (of any index).
In Section 5.2 we provide examples of minimal minitwistor spaces of any index which have a real structure enjoying the conditions in Section 4.
This creates real, positive definite Einstein-Weyl 3-manifolds.  
These minitwistor spaces are obtained as the canonical quotient spaces of the twistor spaces of Joyce's self-dual metrics on the connected sum of complex projective planes.

\vspace{1mm}
\noindent
{\bf Notations and Conventions.}
For a complex space $X$, ${\rm{Sing}}\,X$ means the singular locus of $X$.
For a sheaf $\mathscr F$ on $X$,  we put $h^i(X,\mathscr F)=\dim H^i(X,\mathscr F)$.
By a rational curve, we mean a reduced irreducible curve whose normalization is isomorphic to a complex projective line (as usual).
A nodal curve is a reduced curve which has ordinary nodes as its all singularities.
If $Y$ is a non-singular submanifold in $X$, the normal bundle is denoted by $N_{Y/X}$. (This is used only when $Y\cap {\rm{Sing}}\,X=\emptyset$.)
The base locus of a linear system $|D|$ is denoted by ${\rm{Bs}}\,|D|$.

\section{A Hitchin correspondence and Einstein-Weyl structure on Severi varieties}
\label{Section:Hitchin&Severi}


\subsection{Hitchin correspondence} \label{ss:Hitchin}
First we recall the definition of  Einstein-Weyl structure on complex  manifolds. 
Let $M$ be a complex manifold, and $TM$ and $T^*M$  the holomorphic tangent and cotangent bundles respectively.
A complex metric on $M$ is a holomorphic section 
of the symmetric tensor product $S^2 T^*M$ such that the induced quadratic form on $T_xM$ is non-degenerate for any $x\in M$. 
Two complex metrics $g$ and $g'$ are said to be conformal if there is a non-vanishing holomorphic function $f$ on $M$ satisfying $g'=fg$.
A conformal class of a complex metric $g$ is denoted by $[g]$.
An affine connection on $M$ is a holomorphic connection on $TM$.
An affine connection $\nabla$ on $M$
is said to be compatible with a conformal structure $[g]$ if for each $g\in[g]$, 
there is a (holomorphic) 1-form $a$ on $M$ such that 
\begin{align} \label{compatibility}
 \nabla g = a \otimes g. 
\end{align}
In this case, the pair $([g],\nabla)$ is called a {\em Weyl structure}. 
Further, if $\nabla$ is torsion free, the Weyl structure is called torsion-free. 
On the other hand, a Weyl structure $([g],\nabla)$ is called {\em Einstein-Weyl}\, 
if, for any $g\in[g]$, there is a 
holomorphic function 
$\Lambda$ on $M$ such that 
\begin{align} \label{EWeq}
 R_{(ij)}=\Lambda g_{ij}
\end{align}
holds, where $g_{ij}$ is the metric tensor of $g$ and $R_{(ij)}=\frac{1}{2}(R_{ij}+R_{ji})$
is the symmetrized Ricci tensor of $\nabla$.  

For  {\em real}\,  manifolds,
a Weyl structure and the Einstein-Weyl condition on it is 
defined in a similar way by 
the equations (\ref{compatibility}) and (\ref{EWeq}). 
We say that a Weyl structure $([g],\nabla)$ on a real manifold is positive-definite or negative-definite if $[g]$ is so. 

In this paper, we are concerned with  Einstein-Weyl structures on 3-dimensional manifolds. 
So let  $M$ be a  complex 3-fold 
and  $([g],\nabla)$ a Weyl structure on $M$.
A 2-dimensional subspace $V$ in a tangent space $T_xM$ $(x\in M$) is called a {\em null plane} if 
$[g]$ degenerates on $V$.
The set of null planes in a tangent space $T_xM$ is called the {\em null cone} (at $x\in M$).
A 2-dimensional subspace $V\subset T_xM$ is a null plane iff $V$ tangents to the null cone of $[g]$. 
A 2-dimensional submanifold $\Sigma\subset M$ is called a {\em null surface}\, 
if $T_x\Sigma$ is a null plane for any $x\in \Sigma$.
Then the Einstein-Weyl condition for $([g],\nabla)$ is characterized in terms of null surfaces as follows.

\begin{prop} \label{prop:integrability}(\cite{H82-1}; see also \cite{N08})
 Let $M$ be a 3-dimensional complex manifold and  $([g],\nabla)$
 a torsion-free Weyl structure. 
 Then $([g],\nabla)$ is Einstein-Weyl if and only if, for any $x\in M$ and 
 any null plane $V\subset T_xM$, there is a null surface $\Sigma$ 
 in a small neighborhood of $x$ 
 such that $V$ tangents to $\Sigma$. 
\end{prop}

\noindent Proof. 
Since we later require some detail of the proof, we briefly recall it.
Let $\mathbb{P}(T^*M)\to M$ 
be the  projectivization of the holomorphic cotangent bundle, and consider the null planes bundle
\begin{align} \label{Q(M)}
 Q(M):=\left\{ [\varphi]\in \mathbb{P}(T^*M) \set g(\varphi,\varphi)=0 \right\}, 
\end{align}
where $g\in [g]$, and $\varphi\in T^*M$ is considered to be an element of $TM$ 
by the identification $T^*M\simeq TM$ induced by $g$. 
Then as $[g]$ is non-degenerate, the fiber of the natural projection $\varpi:Q(M)\to M$ is a non-singular conic ($\simeq\mathbb{CP}^1$), 
and each point $[\varphi]\in Q(M)$ defines a null plane 
$\ker\varphi \subset T_{\varpi([\varphi])}M$.  
Since $\nabla$ is compatible with $[g]$, $Q(M)$ is invariant under parallel transports 
of $\nabla$. 
Hence we can define a plane distribution $\mathscr D$ on $Q(M)$ by the condition that   
$\mathscr D_{[\varphi]} \subset T_{[\varphi]}Q(M)$ is the horizontal lift of 
the null plane $\ker\varphi\subset T_{\varpi([\varphi])}M$ with respect to $\nabla$. 
In this situation, we can check by direct calculations that the distribution $\mathscr D$ is integrable 
if and only if $([g],\nabla)$ satisfies the  Einstein-Weyl condition. 
Thus the statement follows. 
Indeed, for a null plane $V=\ker\varphi$, 
we can associate an integral surface
$\tilde\Sigma\subset Q(M)$ of $\mathscr D$ satisfying 
$[\varphi]\in\tilde\Sigma$ iff $\mathscr D$ is integrable. 
Then the image $\Sigma=\varpi(\tilde\Sigma)$ 
is a null surface which tangents to $V$. 
\proofend

\vspace{2mm}
As  explained in the introduction, we shall define an Einstein-Weyl structure on a certain kind 
of complex 3-folds arising from some complex surfaces (which will be called the minitwistor spaces). 
Our method is based on the following construction called 
{\em Hitchin correspondence}  established by Hitchin in \cite{H82-1}. 
Let $S$ be a non-singular complex surface.
A non-singular rational curve $C\subset S$ is called a minitwistor line if $N_{C/S}\simeq \mathscr O_C(2)$ holds. 
As $H^1(N_{C/S})=0$, by Kodaira's theorem, the parameter space $W$ of minitwistor lines becomes a 3-dimensional complex manifold such that for any minitwistor line $D\in W$,  there is a canonical isomorphism $T_DW\simeq H^0(D, N_{D/S})$. 

The complex 3-fold $W$ is naturally equipped with  certain families of 2 and 1-dimensional \ submanifolds as follows. 
First, for any $p\in S$ define $W_p:=\{ D\in W \set p\in D\}$. 
Then by considering the blowing-up $S$ at $p$, it can be shown (\cite{H82-1}) that $W_p$ becomes a 2-dimensional complex submanifold in $W$ if it  is non-empty.  
Next, for any two points $p,q\in S$  we define $W_{p,q}:=\{D\in W \set p,q\in D\}$. 
Then $W_{p,q}$ becomes a 1-dimensional complex submanifold of $W$. 
Note that $W_{p,q}$ naturally makes sense even when $q$ is an infinitely near point of $p$. 
Here for any point $p\in S$, an {\em infinitely near point of $p$} is a point 
on the exceptional curve $E_p$ of the blowing-up of $S$ at $p$. 
When $q$ is an infinitely near point of  $p\in S$, 
then $W_{p,q}$ is defined as $\{ D \in W \set p\in D, q\in  D'\}$ 
where $ D'$ is the strict transform of $D$. 

Then on the parameter space $W$ an Einstein-Weyl structure is defined as follows:

\begin{prop} \label{prop:original} (\cite{H82-1}; see also \cite{N08})
 There is a unique torsion-free Einstein-Weyl structure $([g],\nabla)$ 
 on $W$ satisfying the following properties. 
 (i) The family $\{W_p\}_{p\in S}$ coincides with the set of null surfaces of $[g]$, 
 (ii) The family $\{W_{p,q}\}_{p,q\in S}$ coincides with the set of geodesics. 
 (iii) The family $W_{p,q}$ is a null geodesic iff $p$ and $q$ are infinitely near. 
\end{prop}

\noindent Proof. 
Since we will again require some details of the proof, we give an outline. 
The conformal structure $[g]$ on $W$ is defined in such a way that, for each $D\in W$, 
the null cone 
${\mathscr N}_D\subset T_DW\simeq H^0({\mathscr O}_D(2)) $ is given by 
\begin{align} \label{Hnb}
{\mathscr N}_D=\{\theta\in H^0({\mathscr O}_D(2))\set 
 \text{the equation $\theta=0$ has a double root} \}. 
\end{align}
Then $W_p$ gives a null surface for $[g]$ if it is non-empty. 
Moreover the null surface $W_p$ is totally geodesic for $[g]$. 
Indeed, for each $C\in W_p$, the family $\{W_{p,q}\}_{q\in C}$ gives 
a $\mathbb{CP}^1$-family of geodesics on $W$ 
passing through the point $C$ and contained in $W_p$. 

Let $\varpi:Q(W)\to W$ be the null-plane bundle defined by (\ref{Q(M)}). 
Then we have an isomorphism $Q(W)\simeq \{ (D,p)\in W\times S\set p\in D\}$ 
since each $p\in D$ defines a null plane $T_DW_p\subset T_DW$. 
Hence we obtain the following double fibration: 
\begin{align} \label{cd:Hitchin1}
 \xymatrix{ & Q(W) \ar[dl]_\varpi \ar[dr]^f & \\ 
 W && S, }
\end{align}
where $f$ is the restriction of the projection to the second factor.
From the proof of Proposition \ref{prop:integrability},  null surfaces on $W$ 
naturally lift on $Q(W)$, and $Q(W)$ is foliated by such surfaces. 
Moreover, the leaves of this foliation coincide with  fibers of 
$f$ by construction. 
Note that $W_p=\varpi( f^{-1}(p))$ and that $Q(W)$ is nothing but the 
universal family of minitwistor lines (which are close to $C$). 
Using the facts that each null surface $W_p$ is totally geodesic, 
we can show that there exists a unique torsion-free affine connection $\nabla$ 
such that $\nabla$ is compatible with $[g]$, and 
such that the conditions (ii) and (iii) are satisfied. 
Then $([g],\nabla)$ is Einstein-Weyl by Proposition \ref{prop:integrability}
\proofend


\vspace{2mm}
Finally in this subsection we recall some facts concerning  the null surfaces. 

\begin{rmk} \label{rmk:null_surf} 
{\em (i) If $p$ and $q$ are two points in $S$ which are not infinitely near, 
then $W_{p,q}=W_p \cap W_q$ holds. Hence the intersection of any two 
null surfaces is a non-null geodesic if it is not empty. 
(ii) For each $p\in S$ and $C\in W_p$, there is a unique null geodesic passing through the point
$C$ and contained in $W_p$. Indeed, the tangent line $T_pC\subset T_pS$ 
determines a point $q\in E_p$, where $E_p$ is the exceptional curve of the blowing up of $S$ at $p$. Then $W_{p,q}$ is the required null geodesic.}
\end{rmk}
Notice that (i) and (ii) above are also based on the following basic facts respectively: 
(i)' the intersection of any two null planes at a point is a non-null complex line, and 
(ii)'  each null plane $V$ bijectively corresponds to a null complex line $L$ 
 such that $L\subset V$.

\subsection{Severi varieties} \label{ss:Severi}
First we recall a definition  of Severi varieties and their basic properties. 
The most useful reference on Severi varieties is a book by E.\,Sernesi \cite{SeBook}, especially \S 4.7.
Let $S$ be a non-singular projective algebraic surface and $\mathscr L$ a line bundle over $S$ satisfying $|\mathscr L|\neq\emptyset$.
We often identify the complete linear system $|\mathscr L|$ and  its parameter space $\mathbb PH^0(S,\mathscr L)^*$ (the dual projective space).
For any   integer $\delta>0$, define a subset of $|\mathscr L|$ by
$$
W_{|\mathscr L|,\,\delta}:=
\{C\in |\mathscr L|\set C  {\text{ is reduced, irreducible and Sing\,}}C {\text{ consists of  }} \delta{\text{ ordinary nodes}}\}.
$$
This is called the {\em Severi  variety of $\delta$-nodal curves in $|\mathscr L|$}.
If $C\in  W_{|\mathscr L|,\,\delta}$, we often write $ W_{|C|,\,\delta}$ instead of $  W_{|\mathscr L|,\,\delta}$.
A fundamental result on Severi variety is the following.

\begin{prop}\label{prop-severi}  (i)
The Severi variety $W_{|\mathscr L|,\,\delta}$ is a (possibly empty) locally closed subvariety in $|\mathscr L|$.
(ii) If $C\in W_{|\mathscr L|,\,\delta}$, there is a canonical isomorphism 
\begin{align}
T_CW_{|\mathscr L|,\,\delta}\simeq H^0(\mathscr O_C(C)\otimes\mathscr I_{{\rm{Sing}}\,C}),
\end{align}
where $T_CW_{|\mathscr L|,\,\delta}$ means the Zariski tangent space of $W_{|\mathscr L|,\,\delta}$ at the point $C$, and $\mathscr I_{{\rm{Sing}}\,C}\subset\mathscr O_C$ is the ideal sheaf of \,${\rm{Sing}}\, C$.
(iii) If $H^1(\mathscr O_C(C)\otimes\mathscr I_{{\rm{Sing}}\,C})=0$, 
$W_{|C|,\,\delta}$ is smooth at $C$ and its dimension is equal to $h^0(\mathscr O_C(C)\otimes\mathscr I_{{\rm{Sing}}\,C})$. 
(Namely $H^1(\mathscr O_C(C)\otimes\mathscr I_{{\rm{Sing}}\,C})$ is the obstruction space for deforming $C$ in $S$ preserving all the nodes.)\end{prop}

\begin{rmk}
{\em
(i) The closure of $W_{|\mathscr L|,\,\delta}$ in the projective space $|\mathscr L|$, which becomes an algebraic variety, is also called the Severi variety.
But in this paper we do not adapt it.
(ii) In general,  Severi varieties can become singular, or even non-reduced.
(iii) It is  not necessarily easy to determine non-emptiness of the Severi variety.
Also, even one can show the non-emptiness, it is not easy to determine its connectedness (or irreducibility of the closure).
}
\end{rmk}

By the adjunction formula, the geometric genus of a member $C\in W_{|\mathscr L|,\,\delta}$ is independent of a choice of $C$.
(Namely, it is given by $(1/2)(\mathscr L^2+K_S\cdot \mathscr L)+1-\delta$.)
If it is zero, namely 
when   a Severi variety parametrizes  (nodal) rational curves, it becomes non-singular and has an expected dimension as follows.

\begin{prop}\label{prop-severi2}
Suppose a non-singular projective surface $S$ has a rational curve $C$ which has $\delta\,(> 0)$ nodes as its all singularities.
If $C^2+1-2\delta>0$, the Severi variety $W_{|C|,\,\delta}$  is non-singular and $(C^2+1-2\delta)$-dimensional.
Moreover, under this assumption, $S$ is a rational surface.
\end{prop}

\noindent Proof.
By the Riemann-Roch formula (see \cite[(3.1) Theorem]{BPV}), we have
\begin{align}\label{RR1}
h^0(\mathscr O_C(C)\otimes\mathscr I_{{\rm{Sing}}\,C})-h^1(\mathscr O_C(C)\otimes\mathscr I_{{\rm{Sing}}\,C})=\deg(\mathscr O_C(C)\otimes\mathscr I_{{\rm{Sing}}\,C})+\chi(\mathscr O_C).
\end{align}
Let $\nu:\tilde C\to C$ be the normalization of $C$, so that $\tilde C\simeq\mathbb{CP}^1$.
The nodes of $C$ determine $2\delta$ distinct points on $\tilde C$.
We put $k:=C^2$ for simplicity.
Then as  $\nu^*\mathscr O_C(C)\simeq\mathscr O_{\tilde{C}}(k)$, we have
$$
h^0(\mathscr O_C(C)\otimes\mathscr I_{{\rm{Sing}}\,C})=
h^0(\mathscr O_{\tilde C}(k)\otimes\mathscr O_{\tilde C}(-2\delta))
=k-2\delta+1.
$$
On the other hand, as $C$ has $\delta$ nodes, we have
$$
\deg(\mathscr O_C(C)\otimes\mathscr I_{{\rm{Sing}}\,C})=k-\delta
$$
and
$$
\chi(\mathscr O_C)=h^0(\mathscr O_C)-h^1(\mathscr O_C)=1-\delta.
$$
Hence by \eqref{RR1}, we obtain $H^1(\mathscr O_C(C)\otimes\mathscr I_{{\rm{Sing}}\,C})=0$.
Then by Proposition \ref{prop-severi} (iii) we obtain that the component of $W_{|C|,\delta}$ 
containing $C$ is non-singular and its dimension is given by $h^0(\mathscr O_C(C)\otimes\mathscr I_{{\rm{Sing}}\,C})=C^2+1-2\delta$.

Next we show that $S$ is a rational surface.
If we  $K_S$ denotes  the canonical bundle on $S$, by adjunction formula we have
\begin{align}
2\delta-2=C^2+CK_S.
\end{align}
Hence we have $CK_S=-(C^2+1-2\delta)-1<0$.
Therefore since $C$ actually moves in $S$, the Kodaira dimension of $S$ is $-\infty$.
Hence $S$ is birational to a ruled surface.
If $q:=h^1(\mathscr O_S)=0$, $S$ is a rational surface and we are done.
If $q>0$, 
let $\alpha:S\to T$ be the Albanese map, so that $T$ is a Riemannian surface of genus $q$.
Then obviously a general nodal curve $C$ cannot be contained in a fiber of $\alpha$.
But at the same time, $\alpha|_C$ cannot be surjective since otherwise we obtain a non-trivial map from $\mathbb{CP}^1$ to $T$ via the normalization of $C$.
This means $q=0$.
Hence $S$ is rational.
\proofend

\vspace{2mm}
A more direct and geometric proof of  the first claim of Proposition \ref{prop-severi2} can be given by taking ``the normalization of a tubular neighborhood of $C$\,".
See Section \ref{ss:EWonSeveri} for this.
As an immediate consequence of Proposition \ref{prop-severi2}, we obtain the following  characterization of 3-dimensional Severi varieties of rational curves:

\begin{prop}\label{prop:char60}
If $\,W$ is a 3-dimensional Severi variety of $\delta$-nodal rational curves on a non-singular projective surface $S$,
there exists an integer $m>1$ satisfying $C^2=2m$ and $\delta=m-1$, where $C$ is any one of the curves represented by points of $W$.
Conversely, if a non-singular projective surface $S$ has a $(m-1)$-nodal rational curves $C$ satisfying $C^2=2m$, the Severi variety $W_{|C|,\delta}$ is 3-dimensional.
\end{prop}

As for the linear system $|C|$ on $S$, we have the following

\begin{prop}\label{prop:birat2}
Let $S$ and $C$ be as in Proposition \ref{prop-severi2}, and suppose  $C^2>2\delta-2$.
Then (i) $\dim |C|=C^2+1-\delta$,  (ii)  ${\rm{Bs}}\,|C|=\emptyset$, (iii) the morphism induced by $|C|$ is
birational to the image, whose degree in $\mathbb{CP}^N$ ($N={C^2+1-\delta}$) is equal to $C^2$,
\end{prop}

\noindent Proof.
We again write $k=C^2$. Let $\omega_C$ be the dualizing sheaf.
Then as $C$ is a rational curve with exactly $\delta$ nodes, we have $\deg \omega_C=2\delta-2$.
Hence by the assumption we have $\deg (\omega_C\otimes\mathscr O_C(-C))=2\delta-2-k<0$.
Therefore by duality we have $H^1(\mathscr O_C(C))=0$.
Then by Riemann-Roch formula we obtain 
\begin{align}\label{dim4}
h^0(\mathscr O_C(C))=k+1-\delta.
\end{align}
Hence by the cohomology exact sequence of 
$$
0\lras \mathscr O_S\lras \mathscr O_S(C)\lras \mathscr O_C(C)\lras 0,
$$
and the rationality of $S$, we obtain $h^0(\mathscr O_S(C))=k+2-\delta$, and obtain (i).
Next in order to show (ii) and (iii) we explicitly give a basis of $H^0(\mathscr O_C(C))$.
For this, let $\nu:\tilde C\to C$ be the normalization as before.
Take a non-homogeneous coordinate $z$ on $\tilde C\simeq\mathbb{CP}^1$.
Let $p_1,\cdots,p_{\delta}$ be the nodes of $C$ and $\tilde p_i^1,\tilde p_i^2\in\tilde C$ $(1\le i\le \delta)$ the 2 points determined by the 2 branches at $p_i$.
We write $z=a_i$ and $z=b_i$ for $\tilde p_i^1$ and $\tilde p_i^2$ respectively,
where we can obviously  suppose $a_i\neq\infty$ and $b_i\neq\infty$ for any $i$.
(We do not make distinction for $\tilde p_i^1$ and $\tilde p_i^2$.)
Then by the isomorphism $\nu^*\mathscr O_C(C)\simeq\mathscr O_{\tilde C}(k)$, 
there is an isomorphism $\mathscr O_{\tilde C}(k)_{\tilde p_i^1}\simeq\mathscr O_{\tilde C}(k)_{\tilde p_i^2}$ between 2 fibers of $\mathscr O_{\tilde C}(k)\to\tilde C$.
From these, we obtain an isomorphism
\begin{equation}
H^0(\mathscr O_C(C))\simeq
V:=\{f(z)\in\mathbb C[z]\set \deg f\le k,\, f(a_i)=f(b_i),\,i=1,2,\cdots,\delta\}.
\end{equation}
If we write $f(z)=c_kz^k+c_{k-1}z^{k-1}+\cdots+c_0$, then  the equation $f(a_i)=f(b_i)$ 
gives a homogeneous linear equation for $c_1,c_2,\cdots,c_k$.
For each $\delta<i\le k$, we can choose $f_i\in V$ such that $\deg f_i=i$ (by a dimensional reason).
Then $\{1, f_{\delta+1},f_{\delta+2},\cdots,f_k\}$ are obviously linearly independent.
Hence by \eqref{dim4} these form a basis of $H^0(\mathscr O_C(C))$.
Then since the zero locus of the two sections $1$ and $f_k$ are disjoint, it follows that
the system $|\mathscr O_C(C)|$ is base point free.
Hence, since the restriction map $H^0(\mathscr O_S(C))\to H^0(\mathscr O_C(C))$ is surjective as above, we obtain that $|C|$ is also base point free, meaning (ii).
Let $\phi:S\to\mathbb{CP}^N$ $(N=C^2+1-\delta)$ be the morphism associated to $|C|$.
Then for (iii), since $C\in |C|$, it suffices to show that the morphism $\phi|_C$ is generically $1:1$.
The last morphism is exactly the morphism induced by the system $|\mathscr O_C(C)|$.
Suppose that this morphism is generically $d:1$ with $d>1$.
Let $D\to\phi(C)$ be the normalization of $\phi(C)$.
Then by the universality of the normalization, the composition  $\tilde C\to C\to \phi(C)$ factors as a surjective morphism $\tilde C\to D$ and the normalization $D\to\phi(C)$.
The morphism $\tilde C\to D$ is clearly $d:1$.
Let $w$ be a non-homogeneous coordinate on $D\simeq\mathbb{CP}^1$.
Then the map $\tilde C\to D$ can be written as $w=f(z)$, where $f$ is a polynomial of degree $d$.
Further, the map $D\to\phi(C)\subset \mathbb{CP}^N$ can be written as $w\mapsto (g_1(w),\cdots,g_N(w))$, where $g_i(w)$ is a polynomial of $w$.
Hence the composition $\tilde C\to \phi(C)\subset\mathbb{CP}^N$ is written as
$z\mapsto (g_1(f(z)),g_2(f(z)),\cdots,g_N(f(z)))$, and the degree of the polynomials $g_i(f(z))$ is a multiple of $d$.
Since we have suppose $d>1$, this contradicts our choice of the above basis $\{1, f_{\delta+1},f_{\delta+2},\cdots,f_k\}$.
Hence we have $d=1$, and $\phi$ is a birational morphism.
This implies that the degree of $\phi(S)$ equals to $C^2=k$, as desired.
\proofend

\subsection{Einstein-Weyl structures on the Severi varieties.} \label{ss:EWonSeveri} 

For the purpose of  proving the main result in this section, we first introduce
 notions of `tubular neighborhood' of a nodal curve in a complex surface, and  the `normalization' of the  tubular neighborhood.  
To this end, we first consider a local model of a nodal curve and its normalization. 
Let $(x^1,x^2)$ be the usual coordinate on ${\mathbb C}^2$, and consider 
a nodal curve $Y=Y^1\cup Y^2$ where $Y^l=\{x^l=0\}$. 
Then for small $\varepsilon>0$ and $l=1,2$, we put  $U^l=\{(x^1,x^2)\in\mathbb C^2\set |x^l|< \varepsilon\}$ and call the union  $U^1\cup U^2$ 
as a tubular neighborhood of the nodal curve $Y$.
Then the disjoint union $U^1 \coprod U^2 $ can be regarded as a tubular neighborhood (in the usual sense) of a non-singular non-connected curve $Y^1\coprod Y^2$, and  
 the natural map 
$U^1 \coprod U^2 \to U^1\cup U^2$ gives an extension of 
the normalization $Y^1\coprod Y^2\to Y$.
This is a local model of a tubular neighborhood of nodal curves and  its normalization.

\begin{figure}
\includegraphics{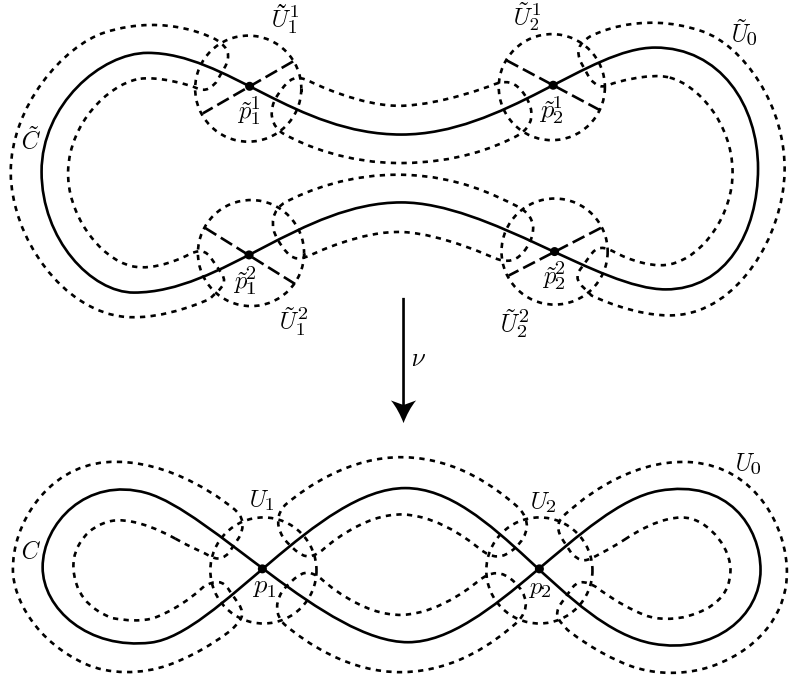}
\caption{the normalization of a tubular neighborhood of a nodal curve}
\label{cover}
\end{figure}

For a global situation,  let $S$ be a compact complex surface, $C\subset S$  a nodal  curve, and 
 $p_1,\cdots,p_n$  the set of nodes of $C$.
Then we can take a covering $\{U_0,U_1,\cdots,U_n\}$ of $C$ by open subsets of $S$ which satisfies the following conditions:
\begin{itemize}
\item[1)] for any $1\le i\le n$ and $0\le j\le n$, $p_i\in U_j$ holds if and only if $j=i$,  
\item[2)] if $0\neq i\neq i'\neq 0$, then $U_i\cap U_{i'}=\emptyset$,
\item[3)] $U_0$ is a tubular neighborhood (in the usual sense) of a non-singular curve $C\backslash(\cup_{i=1}^nU_i)$
\end{itemize}
For any $1\le i\le n$ let $f_i\in \mathscr O(U_i)$ be a defining equation of $C$ in $U_i$.
Then by our choice, $f_i$ can be suppose to be of the form $f_i=x^1_ix^2_i$, where $(x^1_i,x^2_i)$ is a coordinate on $U_i$ such that $p_i=(0,0)$.
Next for any $1\le i\le n$ let $\tilde U_i^1$ and $\tilde U_i^2$ be two copies of  $U_i$, equipped with the same coordinate $(x^1_i,x^2_i)$,   and define $\tilde C_i^1=\{(x^1_i,x^2_i)\in \tilde U_i^1\set x^1_i=0\}$  
and $\tilde C_i^2=\{(x^1_i,x^2_i)\in \tilde U_i^2\set x^2_i=0\}$.
Further let $\tilde p_i^1\in \tilde U_i^1$ and $\tilde p_i^2\in \tilde U_i^2$ be the origins respectively. 
Let $\tilde U_0$ be a copy of $U_0$ and $\tilde C_0\subset \tilde U_0$ be the curve corresponding to $C\cap U_0$.
Then by the above conditions 1) -- 3), the open sets $\{\tilde U_i^1,\tilde U_i^2\set 1\le i\le n\}$ and $\tilde U_0$ naturally glue together 
so that $\{\tilde C_i^1,\tilde C_i^2\set 1\le i\le n\}$ and $\tilde C_0$ also glue together. 
Let $\tilde U$ and $\tilde C$ be the resulting non-singular complex surface and 
the non-singular curve in $\tilde U$ respectively.
Let $$\nu:\tilde U\to U:=\bigcup_{0\le i\le n}U_i$$ be the natural projection.
Then the restriction $\nu|_{\tilde C}$ gives a normalization of $C$.
We call $U$ a tubular neighborhood of $C$, and the map $\nu:\tilde U\to U$ the  normalization of the tubular neighborhood. 
(This construction is illustrated in Figure \ref{cover}). 

The relation between the sheaf $\mathscr O_C(C)\,(=\mathscr O_S(C)|_C)$ and the normal bundle $N_{\tilde C/\tilde U}\simeq \mathscr O_{\tilde C}(\tilde C)
$ is described by the following

\begin{lemma} \label{lem:tang.spaces}
In the above situation, there exists the following exact sequence:
\begin{align} \label{seq:normalization}
0 \lra \mathscr O_{\tilde C}(\tilde C) \lra \nu^*\mathscr O_C(C) 
 \lra \bigoplus_{i=1}^n\left(\mathbb C_{\tilde p_i^1}\oplus\mathbb C_{\tilde p_i^2}\right)
\lra 0, 
\end{align}
where $\mathbb C_p$ means the skyscraper sheaf supported at $p$. 
In particular, there is an isomorphism 
\begin{align} \label{id_tang.spaces}
H^0(\mathscr O_{\tilde C}(\tilde C)) \simeq 
H^0(\mathscr O_C(C)\otimes\mathscr I_{{\rm{Sing}}\,C}). 
\end{align}
\end{lemma}
\noindent Proof. 
We use the above notations prepared for the construction of $\nu:\tilde U\to U$.
We put $\widehat C:=\nu^{-1}(C)$ and $D:=\widehat C-\tilde C$ (subtraction as a divisor).
Then $D$ is a non-compact curve consisting of $2n$ connected components. 
Then if we note an isomorphism $\nu^*\mathscr O_U(C) \simeq \mathscr O_{\tilde U}(\widehat C)$, the sequence  \eqref{seq:normalization} is exactly the third low of the following obvious commutative diagram of exact sequence of sheaves on $\tilde U$:

$$
\begin{CD}
@.0@.0@.0@.\\
@.@VVV@VVV@VVV@.\\0@>>>\mathscr O_{
\tilde U}@>>>\mathscr O_{\tilde U}(D)@>>>\mathscr O_D(D)@>>>0\\
@.@VVV@VVV@VVV@.\\
0@>>>\mathscr O_{
\tilde U}(\tilde C)@>>>\mathscr O_{\tilde U}(\widehat C)@>>>\mathscr O_D(D)\otimes\mathscr O_D\left(\sum_{i=1}^n(\tilde p_i^1+\tilde p_i^2)\right)@>>>0\\
@.@VVV@VVV@VVV@.\\
0@>>>\mathscr O_{
\tilde C}(\tilde C)@>>>\mathscr O_{\tilde C}(\widehat C)@>>>\bigoplus_{i=1}^n\left(\mathbb C_{\tilde p_i^1}\oplus\mathbb C_{\tilde p_i^2}\right)@>>>0\\
@.@VVV@VVV@VVV@.\\
@.0@.0@.0.@.\\
\end{CD}
$$
(In the second low, we have used a relation
$\mathscr O_D(\tilde C)\simeq\mathscr O_D(\sum_{i=1}^n(\tilde p_i^1+\tilde p_i^2))$, which is a consequence of the fact that the two curves $D$ and $\tilde C$ intersect transversally at $\tilde p_i^1$ and $\tilde p_i^2$).

The isomorphism \eqref{id_tang.spaces} is directly deduced from 
the exact sequence \eqref{seq:normalization}. 
\proofend

\vspace{2mm}

Now we are ready to prove the main result in this section.

\begin{thm} \label{thm:definingEWstr}
Any Severi variety $\,W$ of nodal rational curves admits  a natural torsion-free Einstein-Weyl structure, if $\,W$ is 3-dimensional.
\end{thm}

\noindent Proof. 
By Proposition \ref{prop:char60}, there exists an integer $m>1$ such that 
any member $C$ of $W$ is a nodal rational curve in $S$ with $(m-1)$ nodes 
and satisfies $C^2=2m$. 
For each nodal curve $C\in W$, let $\nu:\tilde U_C\to U_C$ be the normalization of a tubular neighborhood of the nodal curve $C$.
As $C^2=2m$ we have $\nu^*\mathscr O_C(C)\simeq\mathscr O_{\tilde C}(2m)$.
Therefore, since the third arrow in the exact sequence \eqref{seq:normalization} is just the evaluation map at $\{\tilde p_i^1,\tilde p_i^2\}$, \eqref{seq:normalization}  means that $\deg \mathscr O_{\tilde C}(\tilde C)= 2m-2(m-1)=2$.

\begin{figure}
\includegraphics{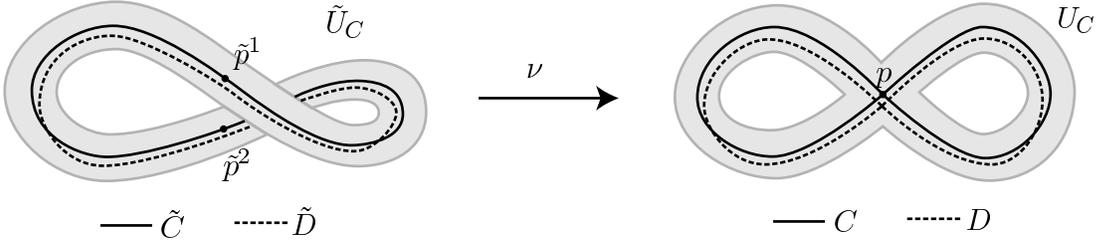}
\caption{the natural lift of the nodal rational curve $D$ near $C$}
\label{tubular2}
\end{figure}

Let $O_C\subset W$ be a sufficiently small open neighborhood of the point $C$ 
such that any $D\in O_C$ is a nodal rational curve contained in $U_C$. 
Notice that each $D\in O_C$ is naturally lifted to a 
non-singular curve $\tilde D$ on $\tilde U_C$ (see Figure \ref{tubular2}). 
Conversely, for each rational curve $\tilde D$ on $\tilde U_C$ which is sufficiently 
close to $\tilde C$, the image $\nu(\tilde D)$ gives a nodal rational curve 
which is a member of $O_C$. 
In this way, $O_C$ coincides with the set of non-singular curves obtained as small deformations of 
$\tilde{C}$ in $\tilde{U}_C$. Hence by the Hitchin correspondence 
(Proposition \ref{prop:original})
 the neighborhood $O_C$ has a natural torsion-free Einstein-Weyl structure.

Now suppose $C_1,C_2\in W$ and $O_{C_1}\cap O_{C_2}\neq\emptyset$. 
We put $O=O_{C_1}\cap O_{C_2}$. 
We claim that the constructed Einstein-Weyl structures on $O_{C_1}$ and 
$O_{C_2}$ agree on $O$. 
Let $\nu_i:\tilde U_{C_i}\to U_{C_i}$ be the normalization of the tubular neighborhood 
$U_{C_i}$ for $i=1,2$. 
Let us fix an arbitrary point $C\in O$. 
We can take a tubular neighborhood $U$ of the nodal curve $C$ such that 
any $D\in O$ is contained in $U$ and  $U\subset U_{C_1}\cap U_{C_2}$. 
Let $\nu: \tilde U\to U$ be the normalization of the tubular neighborhood $U$. 
Then a natural Einstein-Weyl structure on the intersection $O$ is 
induced from $\tilde U$. 
By our construction of $\tilde U\to U$, after shrinking $U$  if necessary, we can define embeddings $\iota_1$ and $\iota_2$ 
so that the diagram 
$$ \xymatrix{
 \tilde U_{C_1} \ar[d]^{\nu_1} \ar@{}[r]|(0.3)\supset & 
 \nu_1^{-1}(U_{C_1}\cap U_{C_2}) \ar[d]^{\nu_1} & 
 \tilde U \ar[d]^\nu  \ar[l]_(0.25){\iota_1}  \ar[r]^(0.25){\iota_2}&
 \nu_2^{-1}(U_{C_1}\cap U_{C_2}) \ar[d]^{\nu_2} &
 \tilde U_{C_2} \ar[d]^{\nu_2} \ar@{}[l]|(0.3)\subset
 \\
 U_{C_1} \ar@{}[r]|(0.45)\supset & U_{C_1}\cap U_{C_2} & U \ar[r] \ar[l] & 
 U_{C_1}\cap U_{C_2} & U_{C_2} \ar@{}[l]|(0.45)\subset
 } $$
commutes. 
Since the Einstein-Weyl structure induced by the Hitchin correspondence is characterized by the local structure of the corresponding complex surface, 
the above diagram indicates that the three surfaces 
$\nu_1^{-1}(U_{C_1}\cap U_{C_2})$, 
$\nu_2^{-1}(U_{C_1}\cap U_{C_2})$ and $\tilde U$ induces the same 
Einstein-Weyl structure on $O$. 
This means that the Einstein-Weyl structures on $O_{C_1}$ and $O_{C_2}$ agree 
on $O$ as required. 
Hence the local structures patch up and define a 
global Einstein-Weyl structure on $W$. 
\proofend

\subsection{The minitwistor spaces.}

In view of Theorem \ref{thm:definingEWstr}, it seems natural to introduce the following

\begin{definition}\label{def:mt}{\em
Let $m>0$ be any integer. Then
by a {\em minitwistor space of index\, $m$}, we mean a pair $(S, |C|)$ of  a non-singular projective algebraic surface $S$ and a complete linear system $|C|$, where $C$ is a nodal rational curve $C$ satisfying $C^2=2m$ which has exactly $(m-1)$ ordinary nodes 
as its all singularities.}
\end{definition}

Evidently,  minitwistor spaces of index 1 are nothing but a (compact) minitwistor spaces in the sense of Hitchin \cite{H82-1}.
In Section 5 we will provide examples of  minitwistor spaces of index $m$ for arbitrary $m>1$.
In general, if $(S,|C|)$ is a minitwistor space of index $m$, there can exist a (nodal rational) curve $C'$ such that $(S,|C'|)$ is a minitwistor space of index $m'\neq m$.
Then the two Severi varieties $W_{|C|,\,m-1}$ and $W_{|C'|,\,m'-1}$ are not biholomorphic in general.
This is a reason why we define the minitwistor space as a pair of $S$ and $|C|$.

If $(S,|C|)$ is a minitwistor space of index $m$, by the results obtained so far, we have the following:
$S$ is  a rational surface by Proposition \ref{prop-severi2}
(since $C^2+1-2\delta=2m+1-2(m-1)=3>0$). 
The Severi variety $W_{|C|,\,m-1}$ ($\subset |C|$) is a 3-dimensional non-singular complex manifold by Propositions \ref{prop-severi2} and \ref{prop:char60}. Furthermore,
it is equipped with a torsion-free Einstein-Weyl structure  by Theorem \ref{thm:definingEWstr}.
We will call any nodal rational curve $C\in W_{|C|,\,m-1}$ as a {\em minitwistor line}. 

Next we show that a blowing up of a minitwistor space is again a minitwistor space.
For this we just need the following

\begin{lemma} \label{lemma:nd1}
Let $(S,|C|)$ be a minitwistor space of (any) index $m$.
Then for any point $p\in S$, there is a member $C\in W_{|C|,\,m-1}$ such that $p\not \in C$.
(Namely, $W\subset\mathbb{CP}^{m+2}$ is non-degenerate.)
\end{lemma}

\noindent Proof.
If $m=1$, this can be readily seen from Kodaira's theorem on displacement of submanifolds.
So let $m>1$.
Suppose that there exists  a point $p\in S$ such that $p\in C$ for any $C\in W_{|C|,\,m-1}$.
Let $\mu:S'\to S$ be the blowing up at $p$.
First assume that general members of $W_{|C|,\,m-1}$ have a node at $p$. 
Let $C'$ be the strict transform of a general member $C\in W_{|C|,\,m-1}$.
Then   $(C')^2=2m-4$  holds.
If $m=2$, $C'$ is a non-singular rational curve satisfying $(C')^2=0$. Hence $|C'|$ is a pencil.
Since we have supposed that $C'$ is a general member (of $W_{|C|,\,m-1}$), this contradicts $\dim W_{|C|,\,1}=3$.
If $m>2$, $C'$ is a rational curve with $\delta':=m-2>0$ nodes.
We have $(C')^2+1-2\delta'=(2m-4)+1-2(m-2)=1>0$. 
Hence by Proposition \ref{prop-severi2}, the Severi variety $W_{|C'|,\,m-2}$ (for $S'$) is 1-dimensional.
This again contradicts $\dim W_{|C|,\,m-1}=3$.
Therefore general members of $W_{|C|,\,m-1}$ are non-singular at $p$.
Hence the strict transform $C'$ of a general member $C\in W_{|C|,\,m-1}$ satisfies $(C')^2=2m-1$.
Then $C'$ is a rational curve with $\delta=m-1$ nodes.
Since $(C')^2+1-2\delta=(2m-1)+1-2(m-1)=2$, Proposition \ref{prop-severi2} again implies 
that $\dim W_{|C'|,\,m-1}=2$. 
This again contradicts $\dim W_{|C|,\,m-1}=3$.
\proofend

\vspace{2mm}

Let $(S,|C|)$ be a minitwistor space of index $m$, and $p\in S$ any point.
Then by Lemma \ref{lemma:nd1} we have $p\not\in C$ for
a general member $C\in W_{|C|,\,m-1}$.
Let $\mu:S'\to S$ be the blowing up at $p$.
Then $C':=\nu^{-1}(C)$ is still an $(m-1)$-nodal rational curve satisfying $(C')^2=2m$.
Further we have $C'\in \mu^*|C|$.
Therefore  the pair $(S',\mu^*|C|)$  becomes a minitwistor space of index $m$.
Let $E$ be the exceptional curve of $\mu$. 
Take any member $C'\in W_{|\mu^*C|,\,m-1}$.  
Then $C'\cap E=\emptyset$ holds since $C'\cdot E=0$ and $C'$ is irreducible by the definition of Severi variety.
This means $\mu(C')\in W_{|C|,\,m-1}$.
Thus we obtain a natural inclusion (as an open subset)
\begin{align}\label{canemb}
W_{|\mu^*C|,\,m-1}\subset W_{|C|,\,m-1}
\end{align}
We remark that in \eqref{canemb} the coincidence does occur in general since $D\in W_{|C|,\,m-1}$ with $p\in D$ is transformed into a reducible curve $D'+E\,(=\mu^*D-E)$ where  $D'$ is the strict transform of $D$, and $D'+E\not\in W_{\mu^*|C|,\,m-1}$ by  the definition of Severi variety.
To exclude those minitwistor spaces (obtained by blowing up of another minitwistor space), we introduce the following

\begin{definition}
\label{def:minimal}{\em
We say that a minitwistor space $(S,|C|)$ is {\em minimal} if it cannot be obtained from another minitwistor space by blowing-up a point; more precisely, if there exist no minitwistor space $(\ol S,|\ol C|)$ and a point $p\in \ol S$ such that $S$ is obtained from $\ol S$ by blowing-up $p$ and such that  $|C|=\mu^*|\ol C|$, where $\mu:S\to \ol S$ is the blow-up at $p$.
}
\end{definition}

Of course, this minimality does not imply that $S$ does not have a $(-1)$-curve.
For minitwistor spaces of index 1, 
the following classification result seems  more or less well known:

\begin{prop}\label{prop:class1}
Let $(S, |C|)$ be a minimal minitwistor space of index 1.
Then  one of the following holds.
(i) $S\simeq\mathbb{CP}^1\times\mathbb{CP}^1$ and $|C|=|\mathscr O(1,1)|$,
(ii) $S\simeq\Sigma_2=\mathbb P(\mathscr O(2)\oplus \mathscr O)$,  and $C$ is a section of the ruling\, $\Sigma_2\to\mathbb{CP}^1$ satisfying $C^2=2$.
\end{prop}

\noindent Proof.
When the index is 1, a standard argument readily implies $\dim |C|=3$ and ${\rm{Bs}}\,|C|=\emptyset$.
Let $\phi:S\to \mathbb{CP}^3$ be the morphism associated to $|C|$. If $|C|$ is composed with a pencil, general member of $|C|$ becomes reducible which cannot happen. Hence $\dim\phi(S)=2$.
Since $C^2=2$ (and $\phi(S)$ cannot be a plane), $\phi$ must be birational over $\phi(S)$ which is a quadratic surface. 
If $\phi(S)$ is a full-rank (i.\,e. non-singular), it is isomorphic to $\mathbb{CP}^1\times\mathbb{CP}^1$, and  $\phi(C)$ is a hyperplane section for any $C\in |C|$.
Hence $\phi(C)\in|\mathscr O(1,1)|$. 
On the other hand since $\phi$ is a birational morphism $\phi$ is decomposed as a composition of the usual blow-ups.
Hence by minimality of $(S,|C|)$, we obtain that $\phi$ is isomorphic. Hence (i) holds.
Alternatively, if $\phi(S)$ is not a full-rank, $\phi(S)$ must be a (quadratic) cone $\ol{\Sigma}_2$ with a unique vertex and $\phi(C)$ is again a hyperplane section.
Hence $\phi$ factors through as  $S\to \Sigma_2\to\ol\Sigma_2$, where $\Sigma_2\to\ol\Sigma_2$ is the minimal resolution of the cone.
But again by minimality, the map $S\to \Sigma_2$ must be isomorphic and  $C$ must be the pull-back of the hyperplane section.
The last curve is 
a $(+2)$-section, which implies (ii).
\proofend

\vspace{2mm}
We note that the proof of (ii) of the proposition means that even if a minitwistor space $(S,|C|)$ is minimal, the morphism $\phi$ associated to $|C|$ can be non-isomorphic (over the image) in general.
If the morphism $\phi$ is non-isomorphic, the minimality means that the (birational) image $\phi(S)$ necessarily has singularities.
(See Remark \ref{rmk-cont} for this.)

In Section 5 we will use the following criterion for the minimality.

\begin{prop}\label{prop:criterion5}
Let $(S,|C|)$ be a minitwistor space of index $m$ and $\phi:S\to\mathbb{CP}^{m+2}$ the birational morphism associated to $|C|$ as before. Then $S$ is minimal if and only if the morphism $\phi$ does not contract any $(-1)$-curves on $S$.
\end{prop}

\noindent Proof.
Suppose that $E$ is a $(-1)$-curve on $S$ such that $\dim\phi(E)=0$.
Let $\mu:S\to\ol S$ be the blowing-down of $E$, and put $\ol C=\mu(C), \,p=\mu(E)$.
Then as ${\rm{Bs}}\,|C|=\emptyset$, $\dim \phi(E)=0$ implies $C\cdot E=0$.
Hence any $C\in W_{|C|,\,m-1}$ is disjoint from $E$.
Therefore $\ol C$ is a nodal rational curve satisfying $\ol C^2=2m$ having $(m-1)$ nodes as its all singularities.
Then $(S,|C|)$ is obtained from $(\ol S,|\ol C|)$ by blowing-up the point $p\in \ol S$.
Hence $(S,|C|)$ is not minimal.
Conversely, let $(S,|C|)$ be non-minimal and $(\ol S, \ol C)$ be a minitwistor space from which $(S, |C|)$ is obtained by blowing-up some point $p\in \ol S$.
Let $E$ be the exceptional curve of $S\to \ol S$.
Then by Lemma \ref{lemma:nd1} we can suppose that $p\not\in\ol C$.
Hence $E\cdot C=0$ holds.
This implies $\dim\phi(E)=0$.
Hence $\phi$ contracts the $(-1)$-curve $E$.
\proofend

\vspace{2mm}

\section{Geometry of Einstein-Weyl structure on the Severi varieties}

\label{Section:geometry}
In the previous section we showed that on 
the 3-dimensional Severi variety  of a minitwistor space (in the sense of Definition \ref{def:mt}) 
there exists a natural Einstein-Weyl structure (Theorem \ref{thm:definingEWstr}). 
In this section, we investigate null surfaces and geodesics on these 3-dimensional complex Einstein-Weyl manifolds.
Throughout this section, $(S,|C|)$ denotes a minitwistor space of (any) index $m$,  $W$ denotes the Severi variety $W_{|C|,\,m-1}$
(which is non-singular and 3-dimensional), and $([g],\nabla)$ denotes the natural Einstein-Weyl structure on $W$.

\subsection{The conformal structure on the Severi varieties} \label{subsect:Tangent}
In this subsection, for any point $C\in W$,
we first identify the tangent space $T_CW$ as a subspace of $H^0(\mathscr O_{\tilde C}(2m))$ and also with $H^0(\mathscr O_{\tilde C}(2))$, where $\tilde C$ is the normalization of $C$.
Next we represent the conformal structure $[g]$ in terms of polynomials on $\tilde C$. 
We fix any $C\in W$.

 First, by Proposition \ref{prop-severi} (i), we have a canonical isomorphism 
\begin{align} \label{can_isom_for_tangent}
T_CW\simeq H^0(\mathscr O_C(C)\otimes\mathscr I_{{\rm{Sing}}\,C}). 
\end{align} 
Let $p_1,\cdots,p_{m-1}$ be the nodes of $C$, and put 
$\nu^{-1}(p_i)=\{\tilde{p}_i^1,\tilde{p}_i^2\}$. 
Next we define 
$\mathscr V_C\subset H^0(\mathscr O_{\tilde{C}}(2m))$ to be the image of 
the composition of the following two canonical injections: 
\begin{align} \label{pullback}
H^0(\mathscr O_C(C)\otimes\mathscr I_{{\rm{Sing}}\,C}) \hookrightarrow 
 H^0(\mathscr O_C(C)) \overset{\nu^*}\hookrightarrow  H^0(\mathscr O_{\tilde{C}}(2m)). 
\end{align}
Then we have an isomorphism $T_CW\simeq \mathscr V_C$. Also we have \begin{align}\label{925}
 \mathscr V_C=\{ s\in H^0(\mathscr O_{\tilde C}(2m))\set 
s(\tilde p_i^1)=s(\tilde p_i^2)=0 \ (1\le i\le m-1)\}. 
\end{align}
Let  $(z_0,z_1)$ be  any homogeneous coordinate on $\tilde{C}\simeq\mathbb{CP}^1$.
 Suppose that in this coordinate the two points $\tilde{p}_i^1$ and $\tilde{p}_i^2$ are represented as
\begin{align} \label{lift_of_node}
\tilde{p}_i^1=(a_i,b_i)\,{\text{ and }}\,\tilde{p}_i^2=(c_i,d_i). 
\end{align}
Of course we have $(a_i,b_i)\neq (c_i,d_i)$ as points on $\mathbb{CP}^1$.
If we put
\begin{align} \label{the_product_f}
f(z_0,z_1):=\prod_{i=1}^{m-1}(b_iz_0-a_iz_1)(d_iz_0-c_iz_1),
\end{align}
then by \eqref{925} each element $s\in \mathscr V_C$ can be written as 
\begin{align}\label{vect1}
 s(z_0,z_1)=f(z_0,z_1)
\left(az_0^2+bz_0z_1+cz_1^2\right).
\end{align}
Here, since  $(a_i,b_i)$ and $(c_i,d_i)$ 
are determined only up to scale, the coefficients $(a,b,c)$ are also determined 
up to scale. 
Thus we obtain  isomorphisms
\begin{align} \label{T_CW1}
T_CW\simeq \mathscr V_C\simeq\{(a,b,c)\set a,b,c\in\mathbb C\}.
\end{align}

Using the normalization of a tubular neighborhood of $C$ explained in the previous section, the tangent space $T_CW$ can be alternatively described as follows.
Let $U$ be a `tubular neighborhood' of the curve $C\in W$, and  let $\tilde{U}\to U$ be  the normalization, which is the extension of the normalization  
$\nu : \tilde{C}\to C$. 
By the proof of Theorem \ref{thm:definingEWstr} and  the Hitchin correspondence explained in Section \ref{ss:Hitchin}, we have an isomorphism
\begin{align} \label{T_CW2}
T_CW\simeq H^0(N_{\tilde{C}/\tilde{U}}) 
\end{align}
where we have $N_{\tilde{C}/\tilde{U}}\simeq\mathscr O_{\tilde{C}}(2)$.
Using the  homogeneous coordinate on $\tilde{C}$ as above, 
each section $\theta\in H^0(N_{\tilde{C}/\tilde{U}})$ 
is written as a quadratic polynomial 
\begin{align} \label{section-theta}
\theta(z_0,z_1)=a z_0^2+b z_0z_1+c z_1^2. 
\end{align}
By Lemma \ref{lem:tang.spaces}, 
the relation between the two isomorphisms $T_CW\simeq\mathscr V_C$ and $T_CW\simeq H^0(N_{\tilde C/\tilde U})$ is clearly given by
\begin{align} \label{id_for_tangents}
 H^0(N_{\tilde{C}/\tilde{U}})\overset\sim\lra \mathscr V_C \ ; \  
 \theta(z_0,z_1) & \longmapsto f(z_0,z_1)\theta(z_0,z_1), 
\end{align}
which is uniquely determined up to scaling.

The null cones of the conformal structure $[g]$ on $W$ can be readily written down using these isomorphisms.
Because $[g]$ was locally defined by regarding $\tilde U$ as 
Hitchin's minitwistor space, the conformal structure $[g]$ is defined, 
as in \eqref{Hnb}, in such a way that 
the null cone 
${\mathscr N}_C\subset T_CW\simeq H^0(N_{\tilde{C}/\tilde{U}})$ 
is given by 
\begin{align} \label{nullcone}
{\mathscr N}_C=\{\theta\in H^0({\mathscr O}_{\tilde{C}}(2))\set 
 \text{the equation $\theta=0$ has a double root} \}. 
\end{align}
By the identification (\ref{id_for_tangents}), 
${\mathscr N}_C$ is also written as
\begin{align} \label{nullcone2}
\mathscr N_C=\{ s\in \mathscr V_C\set \, \,{\text{the quadratic equation }}\,
 s(z_0,z_1)/f(z_0,z_1)=0\,{\text{ has a double root}}\,\}.
\end{align}
If we use the notation of (\ref{T_CW1}), then we can write 
\begin{align}
\mathscr N_C\simeq\{(a,b,c)\in\mathbb C^3\set b^2-4ac=0\}.
\end{align}



\subsection{Null surfaces in the Severi varieties} \label{subsect:Null_surfaces}
In general, Severi varieties are naturally embedded in a projective space by the definition.
If  $W$ is a Severi variety of the minitwistor space $(S,|C|)$, 
$W$ is embedded in $\mathbb{CP}^{m+2}$ since $\dim |C|=m+2$ by Proposition \ref{prop:birat2} (i). 
In this subsection we investigate particular hyperplane sections of $W$, 
which will turn out to be null surfaces of the Einstein-Weyl manifold $(W,[g],\nabla)$. 

First, for arbitrary $p\in S$ we set
$$ |C|_{p}:=\{D\in|C|\set p\in D\} $$ 
which is  a hyperplane in $|C|$ since ${\rm{Bs}}\,|C|=\emptyset$ by Proposition \ref{prop:birat2} (ii), and define a hyperplane section by 
$$ W_p:=W\cap |C|_p. $$
By definition, $W_p=\{D\in W\set p\in D\}$.
By Lemma \ref{lemma:nd1} we have $W_p\neq W$.
However, since $W$ is not closed in the projective space $\mathbb{CP}^{m+2}$,
$W_p$ can be empty.
Actually, we will later show that $W_p=\emptyset$ if there is another point $q\in S$ which satisfies $\phi(p)=\phi(q)$ (Proposition \ref{prop:empty_case}) .




When $W_p$ is non-empty, $W_p$ becomes a divisor on $W$. 
In this case we define a subset of $W_p$ by 
$$W_p^1:=\{D\in W_p\set p\in {\rm{Sing}}\,D\}. $$
Then we have the following

\begin{prop} \label{null_surf}
If $W_p$ is non-empty,  we have the following. 
(i) The complement $W_p\backslash W_p^1$ is a non-singular 
 null surface with respect to $[g]$. 
(ii) The set $W_p^1$ is a non-singular non-null geodesic with respect to $([g],\nabla)$.
(iii) In a neighborhood of each point on $W_p^1$, $W_p$ is a union of two null surfaces (with respect to $[g]$) intersecting transversally along\, $W_p^1$.
 In particular, $W_p$ has ordinary nodes along \,$W_p^1$.  
\end{prop}

\noindent Proof. 
Let $C\in W_p$, and let $\tilde{U}\to U$ be the `normalization of a tubular 
neighborhood' of $C$, that is, the extension of the normalization 
$\nu : \tilde{C} \to C$ to the tubular neighborhoods of $\tilde C$ and $C$. 
By taking a sufficiently small open neighborhood $O\subset W$ of the point $C$, 
 we can assume that any nodal curve $D\in O$ is contained in $U$.  
Recall that each $D\in O$ is naturally lifted to a non-singular curve on $\tilde{U}$ 
which we still denote by $\tilde D$.

Suppose $C\in W_p\backslash W_p^1$ (i.e.\ $p\in C\backslash{\rm{Sing}}\,C$). 
Then $\nu^{-1}(p)$ consists of one point.
Let $\tilde p=\nu^{-1}(p)\in \tilde C$ be the point.
We may assume $O\cap W_p^1=\emptyset$.
By construction, we have 
$$ O\cap W_p=\{ D \in O \set \tilde{p}\in \tilde{D} \}. $$
By Proposition \ref{prop:original} (i), $O\cap W_p$ is a null surface in $O$. 
Hence we obtain (i). 

If $C\in W_p^1$, then $\nu^{-1}(p)$ consists of two points which 
we denote by $\{\tilde{p}^1,\tilde{p}^2\}$. Let us define 
$$\Sigma^k= \{ D\in O \set \tilde{p}^k \in \tilde{D} \},\quad k=1,2.$$
Then again by Proposition \ref{prop:original} (i), $\Sigma^k$ are null surfaces in $O$. Furthermore, by the definition of $W_p$ and $W_p^1$ we have 
\begin{align}\label{branches} O\cap W_p = \Sigma^1 \cup \Sigma^2, \quad 
 O\cap W_p^1= \Sigma^1\cap \Sigma^2.
\end{align} 
Then noting that two different null surfaces always intersect transversally along a non-singular non-null geodesic in general
(Remark \ref{rmk:null_surf} (i)),
\eqref{branches} directly implies (ii) and (iii) of the proposition. 
\proofend

\vspace{2mm}
In the above proof, if $C\in W_p\backslash W_p^1$,  we obtain 
\begin{align} \label{TcWp}
T_CW_p \simeq \{ \theta\in H^0(N_{\tilde C/\tilde U}) \set \theta(\tilde p)=0 \}, 
\end{align}
by the Hitchin correspondence. 
Similarly, if $C\in W_p^1$, then we obtain 
\begin{align} \label{TcSigma}
 T_C\Sigma^k \simeq  \{ \theta \in H^0(N_{\tilde C/\tilde U}) \set 
 \theta(\tilde p^k)=0 \},\quad k=1,2. 
\end{align}
Then any null plane $V$ on $W$ has a unique point $p\in S$ such that $V$ is tangent to the null surface $W_p$ as follows:

\begin{prop} \label{prop:foli}
For any point $C\in W$ and any null plane $V\subset T_CW$, there exists a unique point $p\in S$ which satisfies exactly one of the following: 
(i) $C\in W_p\backslash W_p^1$ and $T_CW_p=V$, or 
(ii) $C\in W_p^1$ and $V$ is the tangent space of one 
 of the two branches of \,$W_p$ at $C$. 
\end{prop}

\noindent Proof. 
Let $\tilde U, U, \tilde C$ be as in the proof of Proposition \ref{prop:foli}. 
By the definition of the conformal structure $[g]$, we can readily see that 
for each null surface $V$ in $T_CW\simeq H^0(N_{\tilde C/\tilde U})$, 
there exists a unique point $\tilde p\in \tilde C$ such that 
$$ V=\{ \theta \in H^0(\mathscr O_{\tilde C}(2m)) \set \theta(\tilde p)=0 \}. $$
We put $p=\nu(\tilde p)$. 
If $p\in C\backslash{\rm Sing}\, C$, then (i) follows from \eqref{TcWp}. 
If $p\in {\rm Sing}\, C$, then the point $\tilde p$ coincides with some 
$\tilde p^k\in \nu^{-1}({\rm Sing}\, C)$ and (ii) follows from \eqref{TcSigma}. 
\proofend

\vspace{2mm}
Next we explain another way to see that $W_p^1$ is exactly the singular locus of $W_p$. 
Take any point $C\in W$. 
We have a natural isomorphism $T_CW\simeq H^0(\mathscr O_C(C)\otimes\mathscr I_{{\rm{Sing}}\,C})$ as in \eqref{can_isom_for_tangent}.
We also have a natural isomorphism $T_C(|C|_p)\simeq H^0(\mathscr O_C(C)\otimes\mathscr I_p)$.
Therefore
at the tangent space level,  we have
\begin{align}\label{523}
T_CW\cap T_C(|C|_p)\simeq H^0(\mathscr O_C(C)\otimes\mathscr I_{{\rm{Sing}}\,C})
\cap
H^0(\mathscr O_C(C)\otimes\mathscr I_p).
\end{align}
If $p\in {\rm{Sing}}\,C$ (namely if $C\in W_p^1$), the right-hand-side is equal to 
 $H^0(\mathscr O_C(C)\otimes\mathscr I_{{\rm{Sing}}\,C})=T_CW$.
This means that the hyperplane $|C|_p$ is tangent to $W$ at the point $C$. 
Hence the hyperplane section $W_p$ is singular at  $C$. 
Therefore $W_p^1\subset{\rm{Sing}}\,W_p$ holds.
On the other hand, if $p\in C\backslash{\rm{Sing}}\,C$ (namely if $C\in W_p\backslash W_p^1$), the right-hand-side of \eqref{523} becomes $
H^0(\mathscr O_C(C)\otimes\mathscr I_{{\rm{Sing}}\,C}\otimes\mathscr I_p)$,
which is readily seen to be 2-dimensional. This means that $|C|_p$ and $W$ intersect transversally at $C$. Therefore $W_p=W\cap |C|_p$ is non-singular at $C$.
Hence we obtain the required coincidence $W_p^1={\rm{Sing}}\,W_p$.

Using the isomorphism \eqref{523}, we can also explain 
that $W_p$ has ordinary nodes along $W_p^1$. 
Fix any $C\in W_p^1$ and let $\nu^{-1}(p)=\{\tilde{p}^1, \tilde{p}^2\}$. 
Now we take $\tilde{q}\in \tilde{C}$ sufficiently close to $\tilde{p}^1$, 
and put $q=\nu(\tilde{q})$. 
Then the surface $W_q$ is non-singular at the point $C$ since 
$q\not\in {\rm{Sing}}\,C$.
Moreover, similarly to (\ref{523}), we have 
\begin{align}
T_CW_q  &\simeq H^0(\mathscr O_C(C)\otimes\mathscr I_{{\rm{Sing}}\,C})
   \cap H^0(\mathscr O_C(C)\otimes\mathscr I_p) \\
  &\simeq H^0 \left(\mathscr O_{\tilde{C}}(2m) \otimes 
  \mathscr I_{\nu^{-1}({\rm{Sing}}\,C)} \otimes\mathscr I_{\tilde q} \right) . 
\end{align}
Hence we obtain 
\begin{align} \label{limit_of_TcWq}
 \lim_{\tilde{q}\to\tilde{p}^1} T_CW_q 
 &\simeq  H^0 \left(\mathscr O_{\tilde{C}}(2m)
   \otimes \mathscr I_{\nu^{-1}({\rm{Sing}}\,C)\backslash \{\tilde p^{1}\} } 
   \otimes (\mathscr I_{\tilde p^1})^2 \right)  \\
&\simeq H^0(\mathscr O_{\tilde C}(2m-(2m-3)-2))
\simeq H^0(\mathscr O_{\tilde C}(1))
\end{align}
which is a 2-dimensional subspace of $\mathscr V_C\simeq T_CW$. 
Obviously, the subspace (\ref{limit_of_TcWq}) coincides with the tangent space $T_C\Sigma^1$, where $\Sigma^1$ is the null surface in the proof of Proposition \ref{null_surf}.
The same argument works for another point $\tilde{p}^2$.



\begin{rmk} \label{rmk:Wp_is_Severi}
{\em
The subvarieties $W_p$ and $W_p^1$ in $W$ are naturally considered as 
Zariski open subsets of certain Severi varieties. 
Indeed, let $\mu : S'\to S$ be the blowing-up at $p$, and  $E_p$  the exceptional curve. 
We fix $C\in W_p\backslash W_p^1$ and let $C'$ be the strict transform of $C$ under $\mu$. 
Then the map 
$$W_p\lra |C'| \ ;\,  D \longmapsto \mu^{-1}(D)-E_p $$ 
gives an (open) embedding of $W_p$ to the 
Severi variety $W_{|C'|,\,m-1}\subset |C'|$
whose dimension is $(2m-1)+1-2(m-1)=2$ by Proposition \ref{prop-severi2}. 
A similar argument works for $W_p^1$. 
When $m=2$, $W_p^1$ is isomorphic to a Zariski open subset of $\mathbb{CP}^1$, since the Severi variety in which $W_p^1$ is embedded is just a pencil. }
\end{rmk}

\subsection{Geodesics on the Severi varieties} \label{subsect:geodesics}
In this subsection, we investigate geodesics on $W$. 
Let $p$ and $q$ be arbitrary two distinct  points on $S$, 
and  set 
$$ \begin{aligned}
  |C|_{p,q} &:= |C|_p\cap |C|_q= \{ D\in |C|\set p,q\in D\}, \\ 
  W_{p,q} &:= W\cap |C|_{p,q}=W_p\cap W_q. \end{aligned} $$
Obviously if $W_p=\emptyset$ or $W_q=\emptyset$, then 
$W_{p,q}=\emptyset$. 
Moreover, we will see later (Proposition \ref{prop:empty_case}) that 
$W_{p,q}$ is empty if $\phi(p)=\phi(q)$, which is the case when 
the coincidence $|C|_{p,q}=|C|_p=|C|_q$ occurs . 
The following proposition means that if $W_{p,q}$ is non-empty, 
$W_{p,q}$ becomes a non-null geodesic on $W$ with respect to $[g]$ 
which has self-intersections. 

\begin{prop} \label{prop:geo1}
Let $p$ and $q$ be two distinct  points in $S$. Suppose $W_{p,q}\neq \emptyset$. 
Let $C\in W_{p,q}$, and let $O\subset W$ be a  neighborhood of the point 
$C$. Then if $O$ is sufficiently small, we have the following. 
 (i) If \,${\rm Sing}\, C \cap \{p,q\}=\emptyset$, then $O\cap W_{p,q}$ is a
 non-null geodesic. 
 (ii) If $p\in{\rm Sing}\, C, q\not\in{\rm Sing}\, C$ or vise versa, 
 then $O\cap W_{p,q}$ is a union of two non-null geodesics intersecting at the point $C$. 
 (iii) If $p,q$ are different nodes of $C$ (so that $m\ge 3$), then $O\cap W_{p,q}$ is a union of four 
 non-null geodesics intersecting at the point $C$. 
\end{prop}

\noindent Proof. 
Let $\tilde{U}\to U$ be the `normalization of a tubular neighborhood' of $C$. 
First suppose $p,q\not\in {\rm Sing}\,C$. 
If we put $\tilde{p}=\nu^{-1}(p)$ and $\tilde{q}=\nu^{-1}(q)$, 
then we have 
$$ O\cap W_{p,q}=\{ D\in O \set \tilde{p},\tilde{q}\in \tilde{D} \}. $$
By Proposition \ref{prop:original} (ii) and (iii), this is a non-null geodesic in $O$. 
Hence (i) holds.

Next suppose $p\in {\rm Sing}\,C$ and $q\not\in {\rm Sing}\,C$. 
Putting $\nu^{-1}(p)=\{\tilde{p}^1,\tilde{p}^2\}$ and $\nu^{-1}(q)=\tilde{q}$, 
we have 
$$ O\cap W_{p,q}=L^1\cup L^2 \qquad \text{where} \quad
 L^j=\{ D\in O \set \tilde{p}^j,\tilde{q}\in\tilde{D}\}. $$
Further $L^1$ and $L^2$ are non-null geodesics by Proposition \ref{prop:original} (iii).
Moreover we have $L^1\cap L^2 = \{C\}$ since $\tilde{C}$ is the unique curve in $O$ 
which contains $\tilde{p}^1,\tilde{p}^2$ and $\tilde{q}\in \tilde{D}$. 
Hence $O\cap W_{p,q}$ is a union of two non-null geodesics intersecting at $C$. 

Finally, suppose $p$ and $q$ are distinct nodes of $C$.  
Putting $\nu^{-1}(p)=\{\tilde p^1,\tilde p^2\}$ and 
$\nu^{-1}(q)=\{\tilde q^1,\tilde q^2\}$, we have 
$$ O\cap W_{p,q}=L^{11}\cup L^{12}\cup L^{21}\cup L^{22} 
\qquad \text{where} \quad
 L^{jk}=\{ D\in O \set \tilde{p}^j,\tilde{q}^k\in\tilde{D}\}. $$
Hence again by Proposition \ref{prop:original} (iii), $O\cap W_{p,q}$ is a union of four non-null geodesics intersecting at $C$. 
\proofend

\vspace{2mm}

Next we show that the hyperplane section $W_p$ can be empty  in general.
\begin{prop} \label{prop:empty_case}
 Let $\phi:S\to \mathbb{CP}^{m+2}$ be the rational map associated to 
 the system $|C|$. 
 If $p$ and $q$ are distinct two points on $S$ satisfying $\phi(p)=\phi(q)$, 
 then the sets $W_p$, $W_q$ and $W_{p,q}$ are all empty. 
\end{prop}
\noindent Proof.
If $\phi(p)=\phi(q)$, then we have $|C|_{p,q}=|C|_p=|C|_q$. 
Hence $W_{p,q}=W_p=W_q$ by definition. 
On the other hand, by Proposition \ref{prop:geo1}, $W_{p,q}$ is 1-dimensional as long as it is non-empty, while $W_p$ and $W_q$ are 2-dimensional if they are non-empty. Therefore since $W_{p,q}\subset W_p$ and $W_{p,q}\subset W_q$, $W_{p,q}=W_p=W_q$ happens only when $W_{p,q}=W_p=W_q=\emptyset$.
\proofend

\vspace{2mm}
Suppose $W_{p,q}$ is non-empty. 
By Proposition \ref{prop:geo1}, a point $C\in W_{p,q}$ is a singular point of $W_{p,q}$ if and only if $({\rm{Sing}}\,C)\cap\{p,q\}\neq\emptyset$. 
We can also prove this in the following way. 
If $({\rm{Sing}}\,C)\cap\{p,q\}=\emptyset$, $C$ is a smooth point of $W_p$ and also of $W_q$ by Proposition \ref{null_surf}. Moreover we have
\begin{align}
T_CW_p\cap T_CW_q&=H^0(\mathscr O_C(C)\otimes\mathscr I_{{\rm{Sing}}\,C}\otimes\mathscr I_p)\cap 
H^0(\mathscr O_C(C)\otimes\mathscr I_{{\rm{Sing}}\,C}\otimes\mathscr I_q)\label{001}\\
&=H^0(\mathscr O_C(C)\otimes\mathscr I_{{\rm{Sing}}\,C}\otimes\mathscr I_{p,q})\label{002}.
\end{align}
The dimension of the last space is $(2m+1)-2(m-1)-2=1$.
This means that the two surfaces $W_p$ and $W_q$ intersect transversally at the point $C$. Hence $W_{p,q}$ is non-singular at the point $C$.
If $({\rm{Sing}}\,C)\cap\{p,q\}\neq\emptyset$, at least one of $W_p$ and $W_q$ has singularities at the point $C$ by Proposition \ref{null_surf}.
Since $W_{p,q}$  is a hyperplane section of $W_p$ and $W_q$, 
$W_{p,q}$ has singularity at $C$.

\vspace{2mm}
So far in this subsection we have supposed $p\neq q$. 
Next we define $W_{p,q}$ when $q$ is an infinitely near point of $p$.
Let $p\in S$, $\mu:S'\to S$ be the blowing-up  at $p$,  $E_p$ the exceptional curve, and $q$  a point on $E_p$.
Then we set 
$$ \begin{aligned}
  |C|_{p,q} &:= \{ D\in |C|\set p\in D, q\in {D'}\}, \\ 
  W_{p,q} &:= W\cap |C|_{p,q}, \end{aligned} $$
where ${D'}$ is the strict transform of $D$. 
By definition we have  $W_{p,q}\subset W_p$. In particular, 
$W_{p,q}$ is empty if $W_p$ is empty. 

\begin{prop} 
Let $p\in S$ and $q\in E_p$ be as above. 
Then if $\,W_{p,q}\neq\emptyset$, 
the set $W_{p,q}$ is a non-singular null geodesic on $W$.  
\end{prop}

\begin{figure}
\includegraphics{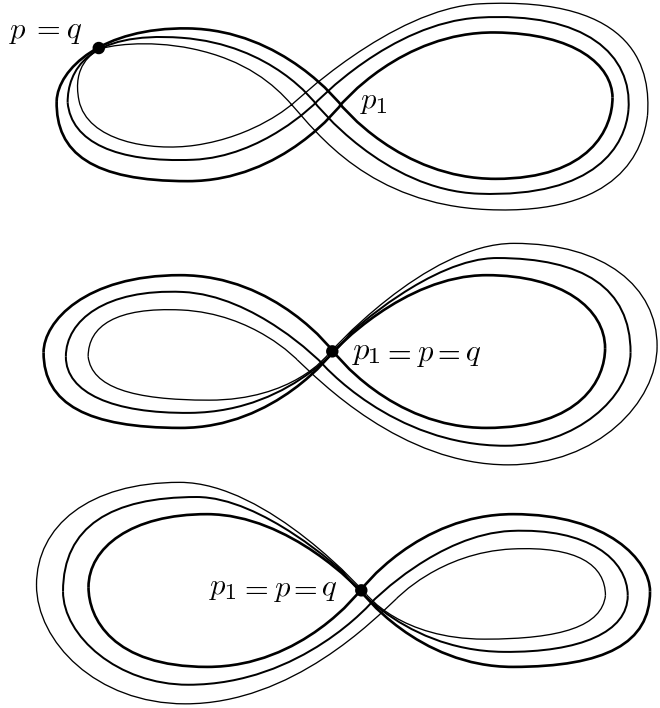}
\caption{Displacements of nodal curves in the case (c)}
\label{fig(c)}
\end{figure}

\noindent Proof. 
Take any $C\in W_{p,q}$. Let $\nu:\tilde C\to C$ be the normalization and 
$\tilde U \to U$  the extension of $\nu$ to the tubular neighborhoods. 
Then we can uniquely define the natural lift $\tilde p\in \tilde C$ and its infinitely near point $\tilde q$ in the following way.
If $p\not\in {\rm Sing}\,C$, then  put $\tilde p=\nu^{-1}(p)$.
Then we define the lift $\tilde q$ to be the point corresponding to the tangent direction $T_{\tilde p}\tilde C$. 
If $p\in {\rm Sing}\,C$, then $\nu^{-1}(p)$ consists of two points. 
We define the lift $\tilde p\in \nu^{-1}(p)$ 
in such a way that $\tilde p$ lies on the branch of which the image by $\nu$ 
tangents to the direction determined by $q$.
The lift $\tilde q$ is defined by the direction $T_{\tilde p}\tilde C$ 
(See Figure \ref{fig(c)}). 
If $O$ is a sufficiently small open neighborhood of the point $C\in W$, 
then by construction, for both cases we have 
$$ O\cap W_{p,q} = 
 \{ D\in O \set \tilde{p}\in\tilde{D}, \tilde{q}\in\tilde{D}' \} $$
where $\tilde{D}'$ is the strict transform of $\tilde{D}$ 
under the blowing-up of $\tilde{U}$ at $\tilde{p}$. 
Hence by Proposition \ref{prop:original} (iii) $W_{p,q}$ is non-singular everywhere, and is a null geodesic. 
\proofend

\begin{rmk}
{\em 
Let $p\in S$ and $q\in E_p$ be as above. 
If $W_p\neq \emptyset$ and  $q\in{\rm{Bs}}\,|\mu^{-1}(D)-E_p|$ for some $D\in W_p$, then we have $|C|_{p,q}=|C|_p$. However this does not happen. 
Indeed, we obtain $W_{p,q}=W_p=\emptyset$ by an argument similar
to the proof of Proposition \ref{prop:empty_case}, and this is a contradiction. 
}
\end{rmk}

\begin{figure}
\includegraphics{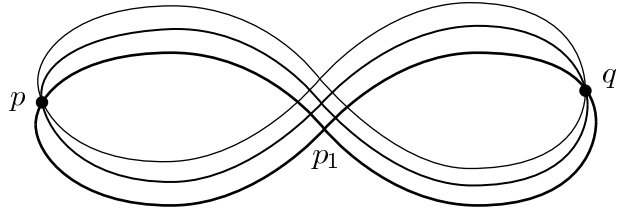}
\caption{Displacement of nodal curves in the case (a1)}
\label{fig(a1)}
\end{figure}

\begin{figure}
\includegraphics{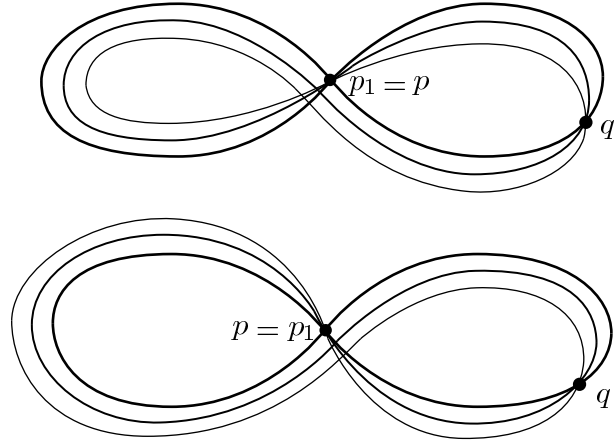}
\caption{Displacements corresponding to the two branches in the case (a2)}
\label{fig(a2)}
\end{figure}

\begin{figure}
\includegraphics{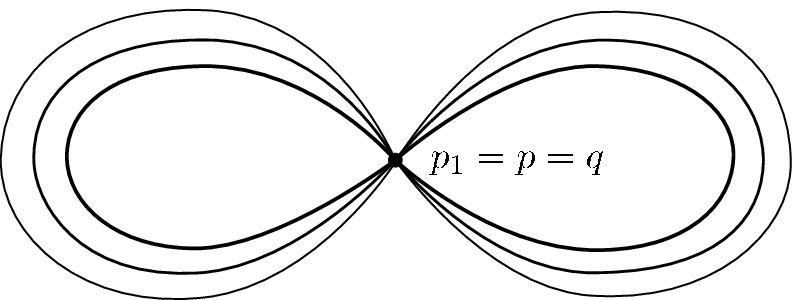}
\caption{Displacement of nodal curves in the case (b)}
\label{fig(b)}
\end{figure}

\vspace{2mm}
Finally in this subsection we illustrate the way how a nodal curve moves when a point on the Severi variety moves along a geodesic.
Take any $C\in W$.
As  in Section \ref{subsect:Tangent}, 
we have $T_CW\simeq H^0(N_{\tilde{C}/\tilde{U}})$. 
Take any tangent direction at $C$, which is represented by a non-zero section $\theta\in H^0(N_{\tilde{C}/\tilde{U}})$. ($\theta$ is uniquely determined up to scaling.)
Let $\tilde{p},\tilde{q}\in\tilde{C}$ be the zeros of $\theta$. 
Putting $\nu(\tilde p)=p$ and $\nu(\tilde q)=q$, the geodesic $\gamma$ which satisfies 
$T_C\gamma=\mathbb C \theta$ is determined (in a neighborhood of $C$) 
in the following way:

\begin{enumerate}
\item[(a1)] If $p\neq q, p\not\in {\rm Sing}\, C$ and $q\not\in {\rm Sing}\, C$, 
  then $\gamma$ coincides with $W_{p,q}$ locally, 
  and $C$ is a smooth point of $W_{p,q}$. (This is the most general case.)
If the point $C$ moves along the geodesic $W_{p,q}$, the nodal  curve $C$ moves as illustrated in Figure \ref{fig(a1)}.

\item[(a2)] If $p\neq q, p\in {\rm Sing}\, C$ and $q\not\in {\rm Sing}\, C$, 
  then $\gamma$ is one of the two branches of $W_{p,q}$. 
If the point $C$ moves along the geodesic $W_{p,q}$, the nodal  curve $C$ moves as illustrated in Figure \ref{fig(a2)}, depending on the two branches of $W_{p,q}$.

\item[(a3)] If $p\neq q, p\in {\rm Sing}\, C$ and $q\in {\rm Sing}\, C$.
  then $\gamma$ is one of the four branches of $W_{p,q}$. 

\item[(b)\phantom{c}] 
  If $p=q$ and $\tilde{p}\neq\tilde{q}$, then
$p\in {\rm Sing}\, C$ and $\tilde p$ and $\tilde q$ correspond to the two branches at $p$.
Hence if $D\in W_{p,q}$, $D$ has a node at $p$.
Therefore $\gamma=W_p^1$.
In this case the nodal curve $C$ moves as in Figure \ref{fig(b)}.

\item[(c)\phantom{b}]
  If $\tilde{p}=\tilde{q}$, then we can consider $q$ as an infinitely near 
  point of $p$, and $\gamma=W_{p,q}$. 
Depending on whether $p\not\in {\rm Sing}\, C$ or $p\in {\rm Sing}\, C$, the nodal curve $C$ moves as illustrated in Figure \ref{fig(c)}.
\end{enumerate}

\begin{rmk} 
{\em
Similarly to Remark \ref{rmk:Wp_is_Severi}, even when $q\in E_p$, 
the geodesics $W_{p,q}$ are naturally considered as Zariski open subsets of some Severi varieties. }
\end{rmk}

\subsection{Double fibration}\label{ss:df}
Finally in this section we explain  a double fibration associated to our construction. 
As our Severi varieties $W$ (of a minitwistor space $(S,|C|)$) carry Einstein-Weyl structure, we have a null plane bundle
 $\varpi : Q(W)\to W$ defined by (\ref{Q(M)}), which is a conic bundle. 
For each point $u\in Q(W)$, let $V_u\subset T_{\varpi(u)}$ be the null plane corresponding to $u$. 
By Proposition \ref{prop:foli}, there exists a unique point $p\in S$ such that
$V_u$ is tangent to $W_p$, where if $\varpi(u)\in W_p^1$, `tangent to $W_p$' means `tangent to one of the branch of $W_p$'. 
Let $f(u):=p$.
Namely, $f(u)\in S$ is the unique point satisfying
$\varpi(u)\in W_{f(u)}$ and 
 $T_{\varpi(u)}W_{f(u)}\supset V_u$.
(When $p\in {\rm{Sing}}\, C$, $T_CW_p$ means the union of two tangent spaces of the two branches of $W_p$ at $C$.)
This way we obtain a map $f:Q(W)\to S$.

Now let us denote $\widetilde W_p:=f^{-1}(p)$ for any $p\in S$.
Then we can obviously write 
$$
 \widetilde{W}_p = \{ u\in Q(W) \set \varpi(u)\in W_p \ \text{and} \ 
 V_u \ \text{is tangent to (one of the branch of)} \ W_p \}.  
$$
This insists that $\widetilde W_p$ is the 
tautological lift of $W_p$. 
In particular, the restriction $\varpi|_{\widetilde W_p}:\widetilde{W}_p\to W_p$ gives the resolution of the singularity $W_p^1$ of $W_p$. 
Clearly $\{\widetilde W_p\set p\in S\}$ foliates $Q(W)$. 
If we denote $Q(W)_C:=\varpi^{-1}(C)$ for the fiber, then we can write as 
$$ Q(W)_C=\{ u\in Q(W) \set C\in W_{f(u)}\} 
  =\{  u\in Q(W) \set f(u)\in C\}. $$
Hence the image $f(Q(W)_C)$ is the nodal rational curve $C$ itself. 
The restriction $f|_{Q(W)_C}:Q(W)_C\to C$  gives the normalization of $C$. 

On the other hand,
since the Severi variety $W$ parametrizes curves on the surface $S$,
there is the universal family which is concretely given by $$R(W)=\{(C,p)\in W\times S \set p\in C \}.$$
For any $C\in W$, the fiber over the point $C$ is exactly the nodal curve $C$.
Obviously the total space $R(W)$ has ordinary nodes along the locus
$\{(C,p)\in R(W)\set p\in {\rm{Sing}}\, C\}$.
Then we can define a natural map $\Psi : Q(W)\to R(W)$ 
by $\Psi(u)=(\varpi(u),f(u))$, which makes the diagram 
\begin{align}\label{diagram4}
 \xymatrix{ & Q(W) \ar[ddl]_\varpi \ar[ddr]^f \ar[d]_\Psi  & \\
  & R(W) \ar[dl] \ar[dr] & \\ 
  W && S. }
\end{align}
commutes. 
The map $\Psi$ gives the normalization of $R(W)$. 
Of course, this diagram corresponds to the diagram \eqref{cd:Hitchin1} in the Hitchin's original case.

\section{Real structure} \label{Section:real}
Throughout this section, as before $(S,|C|)$ always means a minitwistor space of index $m$,  $W=W_{|C|,\, m-1}$ means the (3-dimensional) Severi variety of minitwistor lines, and $([g],\nabla)$ denotes the Einstein-Weyl structure on $W$.
In this section we show that if the rational surface $S$ carries an appropriate real structure, then the component of the real locus of the Severi variety $W$ 
satisfying certain conditions becomes a real, positive definite, 3-dimensional Einstein-Weyl manifolds. 

Let $\sigma$ be
an anti-holomorphic involution on $S$ for which we simply call real structure.  
Obviously, $\sigma$ naturally induces an anti-holomorphic involution on the 
Severi variety $W$ which is also denoted by $\sigma$. 
Set $W^{\sigma}:=\{C\in W\set \sigma(C)=C\}$. 
Since $\sigma$ is an anti-holomorphic involution on the 3-dimensional complex manifold $W$, the real locus $W^{\sigma}$ becomes a real 3-dimensional submanifold (if it is not empty).
Then in general, $W^{\sigma}$ is not necessarily connected, and the type (or the isomorphism class) of the pair $(C,\sigma|_C)$ depends on a choice of the connected component of $W^{\sigma}$ in which $C$ belongs.
For example, even when $m=2$ (i.e. the case of 1-nodal curve), the real locus of $\sigma|_C$ can become a point (i.e. the node) or non-isolated 
(i.e. a 1-dimensional subset which possibly has a node) 
These two curves cannot belong to the same connected component. 
We are interested in the connected components for which the following conditions are satisfied: 

\begin{enumerate}
 \item[$1^{\circ}$] $\sigma(C)=C$, 
fixing each of the $(m-1)$ nodes.
 \item[$2^{\circ}$] 
At any nodes of $C$, $\sigma$ interchanges two branches of $C$.
 \item[$3^{\circ}$]
 $\sigma$ has no real points on $C$ other than the nodes.
 \end{enumerate}
Namely we define 
\begin{align}
 W^\sigma_0:= \{ D\in W^{\sigma}\set 
 D \ \text{satisfies the conditions} \ 1^{\circ},2^{\circ},3^{\circ}\}.
\end{align}
We note that if $C\in W^{\sigma}_0$, then the natural lift of $\sigma|_C$  to  the normalization $\tilde C$ has no real point.
In particular, the real structure induced on $\tilde C$ can be thought as the anti-podal map.

We will prove that the real 3-dimensional manifold $W_0^\sigma$ has a natural positive-definite Einstein-Weyl structure as a {\em real slice} 
of the complex Einstein-Weyl 
structure $([g],\nabla)$ on $W$ defined in Theorem \ref{thm:definingEWstr}, 
where we mean by `real slice' the real objects naturally induced on the real locus by the restriction of complex objects.  
We begin with the conformal structure:

\begin{lemma} \label{lemma:real_conf.str.} The above real 3-manifold \,$W^{\sigma}_0$\! is equipped with a positive-definite conformal structure  which is induced as a real slice of the complex conformal structure $[g]$ on $W$. 
\end{lemma}

\noindent Proof. 
Take any member $C\in W^\sigma_0$.
Recall that we have isomorphism 
$T_CW\simeq \mathscr V_C \subset H^0_{\tilde{C}}(2m)$
(see \eqref{T_CW1}). 
The involution $\sigma$ induces an anti-linear involution on $\mathscr V_C$, and we have $T_C(W^\sigma_0)\simeq \mathscr V_C^\sigma$ where 
$\mathscr V_C^\sigma$ is the fixed subspace under $\sigma$. 
We prove that the restricted conformal structure $[g]|_{\mathscr V_C^\sigma}$ 
naturally defines  a definite real conformal structure. 

Let $(z_0,z_1)$ be a homogeneous coordinate on $\tilde{C}$ 
such that the lifted involution is given by $(z_0,z_1) \mapsto (-\bar z_1,\bar z_0)$. 
We can use a weighted homogeneous coordinate $(z_0,z_1,v)$ on 
$\mathscr O_{\tilde C}(2m)$ with the equivalence relation
\begin{align}
 (z_0,z_1,v) \sim (\lambda z_0, \lambda z_1, \lambda^{2m}v), \quad 
 \lambda \in \mathbb{C^*}, 
\end{align}
where the projection $\mathscr O_{\tilde C}(2m) \to \tilde C$ is given by  
$(z_0,z_1,v) \mapsto (z_0,z_1)$. 
The involution $\sigma$ on $\tilde C$ naturally induces an anti-linear bundle automorphism 
$\tilde{\sigma}:\mathscr O_{\tilde C}(2m)\to\mathscr O_{\tilde C}(2m)$
covering $\sigma$, which can be written as 
\begin{align} \label{inv_of_O(2m)}
\tilde{\sigma}:  (z_0,z_1,v) \longmapsto (-\bar z_1,\bar z_0,h(\bar z_0,\bar z_1)\bar v)
\end{align}
using a holomorphic function $h$ on $\mathbb C^2$. 
We have, however, 
$h(\bar z_0,\bar z_1)=h(\lambda\bar z_0,\lambda\bar z_1)$ 
by the well-definedness. Therefore $h(\bar z_0,\bar z_1)$ is a non-zero constant $h$.
Moreover, since $\tilde\sigma$ is an involution, we obtain $|h|=1$. 

From now on we use the same notation as  
\eqref{lift_of_node}, \eqref{the_product_f} and \eqref{vect1} . 
Take any $s\in \mathscr V_C\simeq H^0(\mathscr O_{\tilde C}(2m))$.
In the above coordinate, the polynomial $s$ defines 
a section of the line bundle $ \mathscr O_{\tilde C}(2m)$ given by
$
  (z_0,z_1) \mapsto  (z_0,z_1, s(z_0,z_1)). 
$
By (\ref{inv_of_O(2m)}), this section is mapped by $\tilde\sigma$ to the 
section represented by the polynomial $s'(z_0,z_1)$ given by 
\begin{align} \label{image_of_tilde_s}
 s'(z_0,z_1) = h\,\overline{s(-\bar z_1,\bar z_0)}
 = h\,\overline{f(-\bar z_1,\bar z_0)}
 \left(\bar a z_1^2 -\bar b z_0 z_1 +\bar c z_0^2)\right. 
{\text{ (by \eqref{vect1})}}.
\end{align}
Now, since the two points $\tilde p_i^1$ and $\tilde p_i^2$ (over the node $p_i\in C$) are antipodal each other, 
we can set $(c_i,d_i)=(-\bar b_i,\bar a_i)$. 
This means 
$$
 \overline{f(-\bar z_1,\bar z_0)}= (-1)^{m-1}\, f(z_0,z_1). 
$$
Substituting this to (\ref{image_of_tilde_s}), and using the isomorphisms 
$T_CW\simeq\mathscr V_C\simeq\{(a,b,c)\set a,b,c\in\mathbb{C} \} $, 
we obtain 
\begin{align}
 \sigma : \mathscr V_C \to \mathscr V_C : (a,b,c) \mapsto 
 (-e^{2i\theta}\bar c, e^{2i\theta}\bar b,-e^{2i\theta}\bar a) 
\end{align}
where we put $(-1)^m h =e^{2i\theta}$. Then if we put 
\begin{align}
 a=\frac{e^{i\theta}}{2}(x_2+ix_3),\ b=e^{i\theta}x_1,\ 
 c=\frac{e^{i\theta}}{2}(-x_2+ix_3),
\end{align}
 we have 
\begin{align}
 T_C(W^\sigma_0)\simeq \mathscr V_C^\sigma \simeq 
 \left\{ (x_1,x_2,x_3) \set x_1,x_2,x_3 \in\mathbb{R} \right\}.
\end{align}
The restriction of the complex conformal structure $[g]$ is represented by 
\begin{align}
 b^2-4ac = e^{2i\theta}\left(x_1^2+x_2^2+x_3^2\right). 
\end{align}
Hence the real slice $[g]|_{\mathscr V_C^\sigma}$ is defined by 
the quadratic form $x_1^2+x_2^2+x_3^2$ which is actually positive-definite. 
\proofend

\begin{lemma}
 The holomorphic affine connection $\nabla$ on $W$ is real on the real locus $W^\sigma_0$ in the sense that 
 $\nabla_XY$ is real for any real  vector fields $X,Y$\! on $W^\sigma_0$. 
\end{lemma}

\noindent Proof. 
For any locally defined holomorphic function $\psi$ on $W$, 
we define another holomorphic function $\tilde\sigma(\psi)$ by 
$\tilde\sigma(\psi)=\overline{\psi\circ\sigma}$. 
Similarly, we define the following involutions:  
\begin{align}
 \tilde\sigma : TW \to TW ;\  X \mapsto \overline{\sigma_*(X)}, 
\end{align}
\begin{align}
 \tilde\sigma : T^*W \to T^*W ;\  a \mapsto \overline{\sigma^*(a)}.
\end{align}
Note that these involutions are fiberwise anti-linear. 
We can check the following:  
\begin{align}
 (\tilde\sigma(X))(\tilde\sigma(\psi))=\tilde\sigma(X\psi), \qquad 
 (\tilde\sigma(a))(\tilde\sigma(X))=\tilde\sigma(a(X)).  
\end{align} 
Now, let us define a connection $\tilde\nabla$ on $TW$ by 
\begin{align} \label{reality_of_nabla}
 \tilde\nabla_XY=\tilde\sigma \left(\nabla_{\tilde\sigma(X)}\tilde\sigma(Y)\right). 
\end{align}
Here, $\tilde\nabla$ is actually a connection, since we have, for example, 
$$ \tilde\nabla_{\psi X}Y
  =\tilde\sigma\left(\nabla_{\tilde\sigma(\psi)\tilde\sigma(X)}
    \tilde\sigma(Y) \right)
  =\tilde\sigma\left(\tilde\sigma(\psi)\nabla_{\tilde\sigma(X)}
    \tilde\sigma(Y) \right)
  =\psi\, \tilde\sigma\left(\nabla_{\tilde\sigma(X)}\tilde\sigma(Y) \right)
  =\psi\tilde\nabla_XY. $$
We claim that:
(i) $\tilde\nabla$ is  torsion-free, (ii) $\tilde\nabla$ is  contained in the projective class $[\nabla]$, and 
(iii) $\tilde\nabla$ is compatible with $[g]$. 
Then since $\nabla$ is the unique connection which satisfies these  conditions
(see Proof of Proposition 4.1 of \cite{N08}),  
we obtain $\tilde\nabla=\nabla$.
Therefore the statement of the lemma follows from the definition (\ref{reality_of_nabla}). 

The claim (i) can be directly checked by using an obvious equality
$[\tilde\sigma(X),\tilde\sigma(Y)]=\tilde\sigma[X,Y]$. 
To check (ii), it is enough to see that every geodesic of $\nabla$ 
is also a geodesic of $\tilde\nabla$. 
(Here, a `geodesic' means an unparametrized geodesic.) 
Take an arbitrary $C\in W$ and a small neighborhood $O\subset W$ of $C$. 
Let $p,q\in C$ be two points which may be infinitely near,  
and  consider the geodesic 
$\gamma :=W_{p,q}\cap O=\{D\in O \set p, q \in D \}$. 
Then $\sigma(\gamma)=W_{\sigma(p), \sigma(q)}\cap \sigma(O)$ 
is also a geodesic. 
If $X$ is a tangent vector field of $\gamma$, 
then $\tilde\sigma(X)$ is a tangent vector field of $\sigma(\gamma)$.
So $\nabla_{\tilde\sigma(X)}\tilde\sigma(X)$ is proportional to $\tilde\sigma(X)$. 
This means that $\tilde\nabla_X X$ is proportional to $X$, which is 
an equivalent condition for $\gamma$ to be a geodesic of $\tilde\nabla$. 

Finally, we check (iii). We fix $g\in[g]$ and define another complex metric $\tilde\sigma(g)$ by 
$(\tilde\sigma(g))(Y,Z)=
 \tilde\sigma\left(g(\tilde\sigma(Y),\tilde\sigma(Z))\right)$. 
Let $a$ be the 1-form satisfying $\nabla g=a\otimes g$, then by a direct calculation we can check that  
\begin{align}
 \tilde\nabla \tilde\sigma(g)=\tilde\sigma(a)\otimes\tilde\sigma(g).
\end{align} 
So $\tilde\nabla$ is compatible with $[\tilde\sigma(g)]$. 
Hence, to prove (iii), it is enough to show the coincidence
$[\tilde\sigma(g)]=[g]$.
For this, notice that each null geodesic is mapped to another null geodesic by $\sigma$. 
Hence for a tangent vector $X$, $g(X,X)=0$ if and only if $\tilde\sigma(g)(X,X)=0$. 
Therefore the null cones of $g$ and those of $\tilde\sigma(g)$ coincide, 
so $\tilde\sigma(g)$ and $g$ are conformally equivalent. 
\proofend

\vspace{2mm}
Now we can show the main result in this section.

\begin{thm} \label{thm:real}
There is a natural torsion-free, positive-definite Einstein-Weyl structure on the real submanifold \,$W^\sigma_0$ \!\!
as a real slice of the complex Einstein-Weyl manifold $(W,[g],\nabla)$. 
\end{thm}

\noindent Proof. 
Take an arbitrary $C\in W^\sigma_0$ and a small neighborhood $O\subset W$ of $C$. 
We can take $g\in[g]$ so that $h=g+\tilde\sigma(g)$ does not vanish on $O$. 
From the proof of the previous lemma, we obtain 
$$[h]=[g], \quad \nabla h = (a+\tilde\sigma(a))\otimes h 
 \qquad \text{on} \ \ O, $$ 
where $a$ is a 1-form satisfying $\nabla g=a\otimes g$. 
Notice that, since $h=\sigma(h)$, $h$ defines a real metric on $W^\sigma_0\cap O$. 
By Lemma \ref{lemma:real_conf.str.}, $h$ is definite
on $W^\sigma_0\cap O$.
By exchanging the sign if necessary,  
we can assume that it is positive-definite. 
Thus we have obtained a torsion-free positive-definite Weyl-structure 
$([h],\nabla)$ on $W^\sigma_0$. 
Moreover, this satisfies the Einstein-Weyl condition since the original $([g],\nabla)$ on $W$ satisfies the condition. 
\proofend

\vspace{2mm} For each $p\in S$ 
we denote $(W^\sigma_0)_p:=W^\sigma_0\cap W_p =\{C\in W^\sigma_0 \set p\in C \}$.
Obviously we have  $(W^\sigma_0)_p=(W^\sigma_0)_{\sigma(p)}$. 

\begin{prop}
For each $p\in S$, the locus $(W^\sigma_0)_p$ is a  geodesic on the real Einstein-Weyl manifold $W^\sigma_0$ if it is not empty.  
\end{prop}

\noindent Proof. 
Take any $C\in W^\sigma_0$ and any complex null plane 
$V\subset T_CW\simeq T_CW^\sigma_0\otimes\mathbb C$. 
Here, we claim $V\neq\sigma(V)$
Let $L\subset V$ be the complex null line 
(see (ii)' written after Remark \ref{rmk:null_surf}). 
Suppose $V=\sigma(V)$, then we obtain $L=\sigma(L)$. 
This means that there exists a real null line $L_{\mathbb R}\subset T_CW^\sigma_0$ 
such that $L=L_{\mathbb R}\otimes \mathbb C$. 
However, such a real null line can not  exist because the conformal structure 
is definite. Hence we obtain $V\neq\sigma(V)$, as required. 
Therefore the intersection $V\cap\sigma(V)$ is a complex line. 
Since $l=V\cap T_CW^\sigma_0$ is the real line satisfying 
$V\cap \sigma(V)=l\otimes \mathbb C$, 
we find that for each null surface $\Sigma\subset W$, 
the real locus $\Sigma\cap W^\sigma_0$ is a real curve if it is not empty. 

Now we put $\bar{p}=\sigma(p)$. 
If $p\neq\bar{p}$, then $p\not\in {\rm Sing}\, C$ by the condition $1^\circ$. 
Hence $(W^\sigma_0)_p$ does not intersects the singular locus  
$W_p^1=\{C \in W \set p\in {\rm Sing}\, C \}$. 
Hence $(W^\sigma_0)_p$ is a real curve. 
Moreover, we have 
$(W^\sigma_0)_p=W_0^\sigma \cap W_p\cap W_{\bar{p}}
 =W^\sigma_0\cap W_{p,\bar{p}}$. 
Since $W_{p,\bar{p}}$ is a geodesic on $W$ and  $\nabla$ is real on 
$W_0^\sigma$, 
$(W^\sigma_0)_p$ is a real geodesic. 

If $p=\bar{p}$, then 
each $C\in W_p^\sigma$ has a node at $p$ by the conditions $1^\circ$ and $3^\circ$. 
Hence $(W^\sigma_0)_p=W^\sigma_0\cap W_p^1$. 
Here, notice that $W_p^1$ is a $\sigma$-invariant complex geodesic in this case.
Hence its real locus $(W^\sigma_0)_p$ is also a geodesic. 
\proofend

\vspace{2mm}
Finally in this section, we explain what happens for the basic diagram \eqref{diagram4} if we take real structure into account.
Let $P(W^\sigma_0):=\mathbb P (TW^\sigma_0)$ be the  projectivization of the tangent bundle of the real manifold $W^{\sigma}_0$. 
Then we can define a natural map given by 
$$ j : Q(W)|_{W^\sigma_0}\longrightarrow P(W^\sigma_0) \,;\ 
 [\varphi] \longmapsto TW^\sigma_0\cap\ker\varphi, $$ 
where $Q(W)|_{W^\sigma_0}$ is the restriction of the null plane bundle $Q(W)$ to the real locus $W^{\sigma}_0$. 
Since $j([\varphi])=j([\bar\varphi])$ and $[\varphi]\neq[\bar\varphi]$ by the conditions $1^{\circ}$--$3^{\circ}$, $j$ is a double cover. 
Combined with the universal family $R(W)\to W$ obtained in Section \ref{ss:df}, we obtain the following commutative diagram: 
\begin{align}\label{diagram5}
 \xymatrix{ P(W^\sigma_0) \ar[dd] 
  & Q(W)|_{W^\sigma_0} \ar[ddl]_\varpi \ar[ddr]^f \ar[d]_\Psi \ar[l]_j
  & \\ 
  & R(W)|_{W^\sigma_0} \ar[dl] \ar[dr] & \\ 
  W^\sigma_0 && \ S. \ }
\end{align}
Recall that 
there is an integrable complex two-plane distribution $\mathscr D$ on 
$Q(W)$ consisting of the fiber directions of $f: Q(W)\to S$. 
Then $\mathscr D$ descends  by $j$ to a real 1-dimensional distribution on $P(W^\sigma_0)$.
The integral curves of the last distribution is the natural lift 
of a geodesic on $W^\sigma_0$. 
Namely, the distribution is the geodesic spray on $W^\sigma_0$.

\section{Explicit examples of the minitwistor spaces}
\label{Section:examples}
In this section, for any $m\ge 2$ we  construct a family of minimal minitwistor spaces of index $m$.
By the results we have obtained so far, the relevant Severi varieties of these minitwistor spaces
have a natural Einstein-Weyl structure.

\subsection{Construction of the minitwistor spaces by blowing-up  $\mathbb{CP}^1\times\mathbb{CP}^1$} 
\label{ss:general_example}
Let $m$ and $k$ be integers satisfying $m\ge 2$ and $1\le k< m$.
Let $D_1$ and $D_2$ be any {\em irreducible} curves on  $\mathbb{CP}^1\times\mathbb{CP}^1$ whose bidegrees are $(k,1)$ and $(m-k,1)$ respectively.
(Of course, these are non-singular rational curves.)
We suppose that the reducible curve $D_1+D_2$ has only ordinary nodes; in other words, we suppose that 
 $D_1$ and $D_2$ intersect transversally at any intersection points.
(So $D_1+D_2$ has exactly $m$ nodes.)
 Next we choose arbitrary $2k$ points $p_1,\cdots,p_{2k}$ on $D_1\backslash D_2$,
 and also $2(m-k)$ points $p_{2k+1},\cdots p_{2m}$ on $D_2\backslash D_1$.
Then let $\mu:S\to \mathbb{CP}^1\times\mathbb{CP}^1$ be the blowing-up at $p_1,\cdots,p_{2m}$.
Let $C_1\subset S$ and $C_2\subset S$ be the strict transforms of $D_1$ and $D_2$ respectively.
We readily have  $C_1^2=C_2^2=0$ and $C_1$ and $C_2$ intersect transversally at $m$ points.
Here, we are allowing (subsets of) the points $\{p_1,\cdots,p_{2k}\}$ or $\{p_{2k+1},\cdots,p_{2m}\}$ to be infinitely near:
in such a case, the blowing-up is always performed at the intersection point of the strict transforms of $D_1$ or $D_2$ with the exceptional curves. 
We show that this surface $S$ satisfies the required property as follows:

\begin{prop}\label{prop-example10}
Let $m\ge 2$ be any integer and
 $S$  (any one of) the rational surface constructed as above. Then there exists a nodal rational curve $C$ on $S$ such that  the pair $(S, |C|)$ becomes a minitwistor space of index $m$ in the sense of Definition \ref{def:mt}.
\end{prop}

It follows from Theorem \ref{thm:definingEWstr} that the Severi variety $W_{|C|,\,m-1}$ (for $S$) is a 3-dimensional complex Einstein-Weyl manifold.

\noindent Proof of Proposition \ref{prop-example10}.
We prove the proposition by showing that for the reducible curve $C_1+C_2$ any one of the $m$ nodes can be smoothed, while all other  nodes remain nodes.
Since $S, C_1$ and $C_2$ are rational satisfying $C_1^2=C_2^2=0$, both of the systems $|C_1|$ and $|C_2|$ are base point free pencils.
Therefore the system $|C_1+C_2|$ is base point free.
By the cohomology exact sequences of  $0\to \mathscr O_S(C_1)\to\mathscr O_S(C_1+C_2)\to\mathscr O_{C_2}(C_1+C_2)\to 0$  and $0\to \mathscr O_S\to\mathscr O_S(C_1)\to\mathscr O_{C_1}(C_1)\to 0$ and noting $\mathscr O_{C_2}(C_1+C_2)\simeq\mathscr O_{C_2}(m)$ and $\mathscr O_{C_1}(C_1)\simeq\mathscr O_{C_1}$, we obtain $H^1(\mathscr O_S(C_1+C_2))=0$.
Further, we easily have $H^2(\mathscr O_S(C_1+C_2))=0$.
On the other hand, by the Riemann-Roch formula we  have
\begin{align}
\chi(\mathscr O_S(C_1+C_2))&=\frac12(C_1+C_2)(C_1+C_2-K_S)+1\\
&=\frac12\left\{(C_1+C_2)^2-C_1K_S-C_2K_S\right\}+1\\
&=\frac12\{2m-(-2)-(-2)\}+1=m+3.
\end{align}
Hence we have $\dim |C_1+C_2|=m+2$.
Let $\phi:S\to\mathbb{CP}^{m+2}$ be the morphism associated to the system $|C_1+C_2|$.
We claim that $\phi$ is birational to its image.
Take different two  points $p$ and $q$ on $S$ which are not on $C_1\cup C_2$.
Then as the pencil $|C_1|$ is base point free, for general $p$ and $q$, there is a curve $C'_1\in |C_1|$ which satisfies $p\in C'_1$ and $q\not\in C'_1$. Then the curve $C'_1+C_2$ satisfies $p\in C'_1+C_2$ and $q\not\in C'_1+C_2$.
This means $\phi(p)\neq\phi(q)$.
Therefore $\phi$ is birational, as claimed.
Then by Bertini's theorem, general members of the system $|C_1+C_2|$ are irreducible (and non-singular).
Further, by adjunction formula, we have $K_SC_1=K_SC_2=-2<0$.
Therefore, any of the nodes of the curve $C_1+C_2$ can be smoothed in the system $|C_1+C_2|$  independently (\cite[page 251]{GLS}).
By taking a smoothing which makes exactly one of the nodes  smooth, we obtain a curve $C\in |C_1+C_2|$ which has exactly $(m-1)$ nodes as its all singularities.
It is obvious that $C$ is irreducible.
Further $C$ must be a rational curve by topological reason.
Finally, as $(C_1+C_2)^2=2m$, we have $C^2=2m$.
Thus we obtain that the curve $C$ satisfies the required properties for the pair $(S,|C|)$ to be a minitwistor space.
\proofend

\vspace{2mm}

The minitwistor space $S$ in Proposition \ref{prop-example10} was obtained as $2m$ points blowing-up of $\mathbb{CP}^1\times\mathbb{CP}^1$.
We note that since $\dim|\mathscr O(k,1)|=2k+1>2k$ (= the number of blowing-up points on $D_1\in|\mathscr O(k,1)|$) and $\dim|\mathscr O(m-k,1)|=2(m-k)+1>2(m-k)$
(= the number of blowing-up points on $D_2\in|\mathscr O(m-k,1)|$), the condition we imposed for the configuration of $2m$ points
(on $\mathbb{CP}^1\times\mathbb{CP}^1$) is an open condition.
Namely, if $p_1,\cdots,p_{2k}$ and $p_{2k+1},\cdots,p_{2m}$ have the irreducible curves $D_1\in |\mathscr O(k,1)|$ and $D_2\in |\mathscr O(m-k,1)|$  on which the points lie respectively, then the same is true for points $p'_1,\cdots,p'_{2m}$ if they are sufficiently close to $p_1,\cdots,p_{2m}$ respectively.
Therefore the above  minitwistor spaces constitute $2\cdot 2m=4m$-dimensional family.
By noting $\dim {\rm{Aut}}(\mathbb{CP}^1\times\mathbb{CP}^1)=6$, 
the number of effective parameters for the family is $4m-6$.

As $m\ge 2$, when the $2m$ points $p_1,\cdots,p_{2m}$ are in a general position, the automorphism group of $S$ is readily seen to be 0-dimensional.
Next by locating the $2m$ points in some special position, 
we provide examples of $S$ which admit an  effective $\mathbb C^*$-action,
or even an effective $\mathbb C^*\times\mathbb C^*$-action. 
For the configuration that allows $\mathbb C^*$-action, we choose  distinct two curves of bidegree $(0,1)$.
Next take $m$ points $p_1,\cdots,p_m$ on one of the $(0,1)$-curves and other $m$ points
$q_{1},\cdots,q_{m}$ on another $(0,1)$-curve.
Here, we are allowing the case that $p_i=p_j$ or $q_i=q_j$ for some $i\neq j$;
in that case, the blowing-up is always performed on the strict transform of the $(0,1)$-curve on which $p_i=p_j$ or $q_i=q_j$ belongs. 
Then the surface $S$ clearly admits an effective $\mathbb C^*$-action which is a lift of the $\mathbb C^*$-action on the first factor of $\mathbb{CP}^1\times\mathbb{CP}^1$.
In order for the surface $S$ to have the required nodal rational curves,  we suppose that the $2m$ points satisfy the following (genericity) condition:

\vspace{1.5mm}\noindent
($\ast$) Let $a_i:=\pi_1(p_i)\in\mathbb{CP}^1$ and $b_i:=\pi_1(q_i)\in\mathbb{CP}^1$, where $\pi_1:\mathbb{CP}^1\times\mathbb{CP}^1\to\mathbb{CP}^1$ denotes the projection to the first factor.
Then after possible renumbering for the indices of $a_i$ and $b_i$ independently, there exists a subset $I\subset \{1,2,\cdots,m\}$ with $I\neq \emptyset$ and $I\neq\{1,\cdots,m\}$ 
such that $\{a_i\set i\in I\}\cap \{a_i\set i\not\in I\}=
\{b_i\set i\in I\}\cap \{b_i\set i\not\in I\}=
\{a_i\set i\in I\}\cap \{b_i\set i\in I\}=\{a_i\set i\not\in I\}\cap \{b_i\set i\not\in I\}=\emptyset$.
(See Remark \ref{rmk:exc} for an example of the configuration which does not satisfy this condition.) 

\vspace{1.5mm}
\noindent Note that if the $2m$ points $\{p_i,q_i\set 1\le i\le m\}$ satisfies this condition, then 
$\{p'_i,q'_i\set1\le i\le m\}$ also satisfies the condition if $p'_i$ and $q'_i$ are sufficiently close to $p_i$ and $q_i$ respectively, and if $p'_i$ and $q'_i$ are lying on the same $(0,1)$-curves respectively.

\begin{prop}\label{prop-C*}
Assume that  the $2m$ points on $\mathbb{CP}^1\times\mathbb{CP}^1$ satisfy the above condition  $(\ast)$. Then the surface $S$ (with $\mathbb C^*$-action)  has a rational curve $C$ such that the pair $(S,|C|)$ is a minimal minitwistor space of index $m$.
\end{prop}

Since the $\mathbb C^*$-action on $S$ clearly induces a non-trivial $\mathbb C^*$-action on the set of nodal curves, the Severi variety $W_{|C|,m-1}$ on  $S$ has a non-trivial $\mathbb C^*$-action. Namely,   $W_{|C|,\,m-1}$ is an Einstein-Weyl manifold which admits a $\mathbb C^*$-action.
The last $\mathbb C^*$-action preserves the Einstein-Weyl structure, since the action clearly preserves both of the set of null cones \eqref{nullcone} (or \eqref{nullcone2}) and the set of  geodesics $W_{p,q}$.

\vspace{1.5mm}
\noindent Proof of Proposition \ref{prop-C*}.
As in the proof of Proposition \ref{prop-example10}, it suffices to find 2 non-singular rational curves $C_1$ and $C_2$ on $S$ satisfying $C_1^2=C_2^2=0$ and $C_1C_2=m$, and intersecting transversally.
We give these curves explicitly.
Let $u$ and $v$ be non-homogeneous coordinates on the first and second factor of $\mathbb{CP}^1\times\mathbb{CP}^1$ respectively.
In these coordinates we can suppose that the $2m$ points are explicitly defined by $(u,v)=(a_i,0)$ for $p_i$ ($1\le i\le m$) and $(u,v)=(b_i,\infty)$ for $q_i$  ($1\le i\le m$).
Then  we  define curves $D_1$ and   $D_2$ by the equations
\begin{align}\label{defcurves}
v=\frac{\prod_{i\in I}(u-a_i)}{\prod_{i\in I}(u-b_i)}\hspace{2mm}{\text{ and }}\hspace{2mm}
v=c\,\frac{\prod_{i\not\in I}(u-a_i)}{\prod_{i\not\in I}(u-b_i)}
\end{align}
respectively, where $c\in\mathbb C^*$.
By the condition $(\ast)$ we have  $\{a_i\set i\in I\}\cap\{b_i\set i\in I\}=\emptyset$
and  $\{a_i\set i\not\in I\}\cap\{b_i\set i\not\in I\}=\emptyset$.
Therefore the denominators and numerators in \ref{defcurves} have no common factor.
Hence if  $k$ means the number of elements of $I$, so that $1\le k<m$, the bidegrees of $D_1$ and $D_2$ are $(k,1)$ and $(m-k,1)$ respectively.
Clearly, these are  irreducible curves, $\{p_i,q_i\set i\in I\}\subset D_1$,
and  $\{p_i,q_i\set i\not\in I\}\subset D_2$.
By the condition $(\ast)$, $\{p_i,q_i\set i\in I\}\cap \{p_i,q_i\set i\not\in I\}=\emptyset$ holds.
In particular, $D_1$ and $D_2$ do not intersect on the two $(1,0)$-curves on which $p_i$ and $q_i$ belongs.
Hence, for general choices of  $c\in \mathbb C^*$, $D_1$ and $D_2$ intersect transversally (at $m$ points).
Then as the two curves $C_1$ and $C_2$, it suffices to  choose the strict transforms of $D_1$ and $D_2$ respectively.
\proofend

\begin{rmk}\label{rmk-cont}{\em
For the surface $S$ and the curve $C$ in Proposition \ref{prop-C*}, the linear system $|C|$ is base point free and $(m+2)$-dimensional by Proposition \ref{prop:birat2}.
If $\phi:S\to\mathbb{CP}^{m+2}$ still denotes the morphism associated to $|C|$, $\phi$ contracts two rational curves which are the strict transforms of the $(0,1)$-curves on which the $2m$ points lie.
The image of these two curves become cyclic quotient singularities of $\phi(S)$.
}\end{rmk}

\begin{rmk}\label{rmk:exc}{\em
A simple example of a configuration of the $m+m$ points (on $\mathbb{CP}^1\times\mathbb{CP}^1$) which does not satisfy the condition $(\ast)$ is the following.
Let $m=3$, and let the $2m=6$ points to be  $\{(0,0), (1,0),(\infty,0),(0,\infty),(1,\infty),(\infty,\infty)\}$.
For this configuration, the proof of Proposition \ref{prop-C*} does not work.
We do not know whether for this kind of configuration the surface $S$ has the required nodal rational curves or not.
}\end{rmk}

For general choices of the $2m$ points satisfying the condition $(\ast)$, the identity component of the automorphism group of the minitwistor space $S$ is clearly $\mathbb C^*$.
By a similar consideration for the case of general configurations (without $\mathbb C^*$-action),
the number of effective parameters for the family of the minitwistor spaces in Proposition \ref{prop-C*}
is given by $2m-3$.

Among all configurations of the $2m$ points satisfying $(\ast)$, the surface $S$ becomes a toric surface iff the set $\{a_i,b_i\set 1\le i\le m\}\,(\subset\mathbb{CP}^1)$ consists of exactly $2$ points;
in other words, iff $p_i=(0,0)$ for $1\le i\le k$, $p_i=(\infty,0)$ for $k< i\le m$, $q_i=(0,\infty)$ for $k< i\le m$ and $q_i=(\infty,\infty)$ for $1\le i\le k$
(after renumbering the indices and changing the coordinates).
By looking the self-intersection numbers of the irreducible components of the unique $\mathbb C^*\times\mathbb C^*$-invariant anticanonical curves (which determines the toric surface uniquely), it can be readily verified that the structure of these toric surfaces is independent of $k$.
Therefore, for each $m\ge 2$ we have obtained exactly one toric surface whose Severi variety of $(m-1)$-nodal rational curves is a 3-dimensional complex Einstein-Weyl manifold.
Exactly as in the case with $\mathbb C^*$-action, $\mathbb C^*\times\mathbb C^*$ acts on these manifolds preserving the Einstein-Weyl structure.
Therefore for each $m\ge 2$ we obtain a  complex Einstein-Weyl 3-fold  which admits a $\mathbb C^*\times\mathbb C^*$-action.

Next we show the minimality of the minitwistor spaces obtained so far:
\begin{prop}
All the minitwistor spaces $(S,|C|)$ in Propositions \ref{prop-example10} and \ref{prop-C*} are minimal in the sense of Definition \ref{def:minimal}.
\end{prop}
\noindent Proof. Let $\mu:S\to\mathbb{CP}^1\times\mathbb{CP}^1$ be the blowing-up,  $E_1,\cdots,E_{2m}$ the exceptional curves which satisfy $E_i\cdot E_j=-\delta_{ij}$ for $1\le i,j\le 2m$, and $\phi:S\to\mathbb{CP}^{m+2}$ the birational morphism induced by $|C|$ as in the  proof of Proposition \ref{prop-example10}.
(In particular, $E_i$'s are not necessarily irreducible.)
For the minimality, by Proposition \ref{prop:criterion5} it suffices to show that $\phi$ does not contract any $(-1)$-curves on $S$.
Suppose that $E$ is a $(-1)$-curve  which is contracted by $\phi$.
If $E=E_j$ for some $1\le j\le 2m$, $E_j$ is irreducible, and we have 
\begin{align}
C\cdot E=\left(\mu^*\mathscr O(m,2)-\sum_{i=1}^{2m}E_i\right)\cdot E_j=-E_j^2=1.
\end{align}
Hence (together with the fact that  ${\rm{Bs}}|C|=\emptyset$)  we have $\dim\phi(E)=1$.
Hence $E\neq E_j$ for any $1\le j\le 2m$.
So we can write
\begin{align}
E\sim
\mu^*\mathscr O(k,l)-\sum_{i=1}^{2m}m_iE_i,
\end{align}
where $\sim$ denotes linearly equivalence, $k,l$ and $m_i$ satisfy $k\ge0, \,l\ge 0, \,k+l>0, m_i\ge 0$ and $ \sum_{i=1}^{2m} m_i>0$.
Then as $E$ is a $(-1)$-curve, we have 
\begin{align}
-1=K_S\cdot E&=\left(\mu^*\mathscr O(-2,-2)+\sum_{i=1}^{2m} E_i\right)\cdot \left(\mu^*\mathscr O(k,l)-\sum_{i=1}^{2m}m_iE_i\right)\notag\\
&=-2k-2l+\sum_{i=1}^{2m}m_i.\label{147}
\end{align}
On the other hand, as we have supposed $\dim\phi(E)=0$, we have 
\begin{align}
0=C\cdot E&=\left(\mu^*\mathscr O(m,2)-\sum_{i=1}^{2m} E_i\right)\cdot \left(\mu^*\mathscr O(k,l)-\sum_{i=1}^{2m}m_iE_i\right)\notag\\
&=ml+2k-\sum_{i=1}^{2m}m_i.\label{148}
\end{align}
By \eqref{147} and \eqref{148} we obtain $l(m-2)=-1$.
As $l$ and $m-2$ are non-negative integers, this means $l=m=1$.
This contradicts our assumption $m>1$.
Therefore there exists no $(-1)$-curve $E$ on $S$ which is contracted by $\phi$.
\proofend

\subsection{Examples of the minitwistor spaces with a real structure and real twistor lines} \label{ss:real_example}
The complex surfaces we shall consider next are the rational surfaces given in \cite[Section 2]{Hon-mt}.
Though they can be shown to be included in the examples in the last subsection, 
an advantage of the  surfaces in this subsection is that they are equipped with a natural real structure which is induced from that of the twistor spaces (of  real self-dual 4-manifolds), and we can find real nodal rational curves satisfying the conditions we have considered throughout Section 4.
By the result in Section 4, this means that the real locus of the  Severi varieties become  real (3-dimensional) Einstein-Weyl manifolds.

First we briefly recall some of the results in \cite{Hon-mt}.
We consider (arbitrary) effective $U(1)^2$-action on $n\mathbb{CP}^2$ and take any one of Joyce's self-dual metrics on $n\mathbb{CP}^2$ which are invariant under the $U(1)^2$-action \cite{J95}.
Let $Z$ be the twistor space of the self-dual metric.
The $U(1)^2$-action on $n\mathbb{CP}^2$ naturally induces a holomorphic $G:=\mathbb C^*\times\mathbb C^*$-action on $Z$.
On the other hand, the  $U(1)^2$-action on $n\mathbb{CP}^2$ has exactly $(n+2)$ invariant two-spheres.
As the next step for obtaining the required complex surface, 
choose any one of these $U(1)^2$-invariant spheres and
let $K_1\subset U(1)^2$ be the isotropy subgroup which fixes any points on the  sphere.
 $K_1$ is isomorphic to $ U(1)$.
Let $G_1\subset G$ be the complexification of $K_1$.
We have $G_1\simeq\mathbb C^*$.

Let $F$ be the canonical square root of the anticanonical line bundle of $Z$.
Then the $G$-action on $Z$ naturally lifts on $F$, so that also on the tensor product $kF$, $k>0$.
Let $H^0(Z, kF)^G$ (resp. $H^0(Z, kF)^{G_1}$) be the subspace consisting of all $G$-invariant (resp. $G_1$-invariant) sections, 
and $|kF|^G$ (resp. $|kF|^{G_1}$) the corresponding linear system.
$|F|^G$ is a pencil whose members are smooth toric surfaces.

Under this situation, we have the following.

\begin{prop}\label{prop:m}(\cite[Section 2]{Hon-mt})
There exists a unique integer $m$ satisfying the following.
(i) If $k<m$, the system $|kF|^{G_1}$ is composed with the pencil $|F|^G$, so that $\dim |kF|^{G_1}=k$.
(ii) $\dim |mF|^{G_1}=m+2$.
(iii) If $\Phi_m^{G_1}:Z\to\mathbb{CP}^{m+2}$ denotes the rational map associated to the system $|mF|^{G_1}$, the image $\Phi_m^{G_1}(Z)$ is a normal rational surface whose degree (in $\mathbb{CP}^{m+2}$) is $2m$.
\end{prop}

Note that the integer $m$ is explicitly computable through the algorithm given in \cite[Section 2, Procedure (A)]{Hon-mt}.
As in \cite[Def.\,2.9]{Hon-mt} we write $\mathscr T:=\Phi_m^{G_1}(Z)$.
Recall that $\mathscr T$  can be regarded as a quotient space with respect to the $G_1$-action: namely for general points of $\mathscr T$, the inverse image under $\Phi_m^{G_1}$ are the closure of the $G_1$-orbits.
(More invariantly, $\mathscr T$ is exactly the `canonical quotient space' of the Moishezon twistor space $Z$ under the $G_1$-action, as proved in \cite[Appendix]{Hon-mt}).
Also recall that $\mathscr T$ has a natural real structure induced by that on $Z$.
Although $\mathscr T$ always has isolated singularities as long as $m\ge 2$ (\cite[Prop.\,2.14]{Hon-mt}), we have the following.

\begin{prop}\label{prop-existence}
On the rational surface $\mathscr T$ there exists a curve $\mathscr C$ satisfying the following.
(i) $\mathscr C\cap{\rm{Sing}}\, \mathscr T=\emptyset$.
(ii) $\mathscr C$ is a rational curve having $(m-1)$-nodes as its all singularities.
(iii) $\mathscr C^2=2m$. 
(iv) $\mathscr C$ is real, satisfying the conditions $1^{\circ}$, $2^{\circ}$ and $3^{\circ}$ in Section 4.
\end{prop}

\noindent Proof.
First we  recall the structure of $Z$ more closely.
By \cite[Prop.\,2.11]{Hon-mt}, the system $|mF|^G$ is composed with the pencil $|F|^G$.
As in \cite[Diagram (13)]{Hon-mt} we  have the following commutative diagram of meromorphic maps:
\begin{equation}\label{cd1}
 \CD
Z@>{\Phi_m^{G_1}}>>\mathbb P(H^0(mF)^{G_1})^*\\
 @V\Psi_m VV @VV{\pi_m}V\\
\Lambda_m@>{\iota}>>\mathbb P(H^0(mF)^G)^*,\\
 \endCD
 \end{equation}
where $\Psi_m$ is the map associated to the system $|mF|^G$, $\Lambda_m$ is the image of $\Psi_m$ so that $\Lambda_m\simeq\mathbb{CP}^1$, and 
$\pi_m$ is the linear projection associated to the obvious inclusion $H^0(mF)^G\subset H^0(mF)^{G_1}$
whose fibers are $\mathbb{CP}^2$.
Here, $\iota$ embeds $\Lambda_m\simeq\mathbb{CP}^1$ as a rational normal curve,
and fibers of the restriction $\pi_m|_{\mathscr T}:\mathscr T\to\iota(\Lambda_m)$ are conics, by which $\mathscr T$ has a structure of (rational) conic bundle  over $\iota(\Lambda_m)\simeq\mathbb{CP}^1$. 

To find the curve $\mathscr C$ in the proposition, take a real twistor line $L$ which is disjoint from ${\rm{Bs}}\,|F|^G$.
(Recall that ${\rm{Bs}}\,|F|^G$ is exactly the cycle of rational curves which is the unique $G$-invariant anticanonical curve on a smooth member of $|F|^G$.)
We will show that if $L$ is  sufficiently general then the image $\mathscr C:=\Phi_m^{G_1}(L)$ satisfies the required properties in the proposition.

For verifying the  property (i), recall that $\rm{Sing}\,\mathscr T$ consists of (a) one conjugate pair $\{P_{\infty},\ol{P}_{\infty}\}$ of cyclic quotient singularities and (b) real singularities (\cite[Prop.\,2.14]{Hon-mt}).
For the former ones, we have $(\Phi_m^{G_1})^{-1}(\Phi_m^{G_1}(P_{\infty}))\subset {\rm{Bs}}\,|F|^G$, and hence also $(\Phi_m^{G_1})^{-1}(\Phi_m^{G_1}(\ol P_{\infty}))\subset {\rm{Bs}}\,|F|^G$.
For the latter ones, the inverse images (under $\Phi_m^{G_1}$) of the singularity is  one of the $G$-invariant twistor lines, so that it intersects  ${\rm{Bs}}\,|F|^G$.
Since $L\cap {\rm{Bs}}\,|F|^G=\emptyset$ by our assumption, these imply  $\mathscr C\cap{\rm{Sing}}\, \mathscr T=\emptyset$.
Hence $\mathscr C$ satisfies (i).

Next to show (iii) we claim that $\mathscr C$ is contained in a hyperplane in $\mathbb P(H^0(mF)^{G_1})^*$.
This can be proved in the same way as in \cite[Prop.\,4.4 (d)]{Hon-mt}.
Actually as $F\cdot L=2$, $\Psi_m|_L$ is $2:1$ over $\Lambda_m$. Then by the commutativity of the diagram \eqref{cd1}, $(\pi_m\circ\Phi_m^{G_1})|_{L}$ is  $2:1$ over $\iota(\Lambda_m)$.
Hence $\pi_m|_{\mathscr C}$ is either generically $1:1$ or $2:1$.
But since $\mathscr C$ intersects any irreducible component of reducible fibers of $\pi_m|_{\mathscr T}$ at least once as $L$ intersects any irreducible component of reducible members of $|F|^{G}$ exactly once, and the last irreducible component is mapped to an irreducible component of a reducible fiber of $\pi_m|_{\mathscr T}$,
we deduce $\pi_m|_{\mathscr C}$ is  $2:1$ over $\iota({\Lambda}_m)$.
Furthermore, by the property (i), we have $\mathscr C\cap{\rm{Sing}}\,\mathscr T=\emptyset$.
On the other hand, as $\pi_m|_{\mathscr T}:\mathscr T\to\Lambda_m$ is a conic bundle, we see that hyperplane sections of $\mathscr T$ also have the same intersection numbers with irreducible components of any fibers of $\pi_m|_{\mathscr T}$.
Moreover, general hyperplane sections do not intersect ${\rm{Sing}}\,\mathscr T$.
These imply that the curve $\mathscr C$ and  hyperplane sections of $\mathscr T$ determine the same element in $H^2(\mathscr T,\mathbb Z)$.
Hence $\mathscr C$ is contained in a hyperplane section, as claimed.
Therefore we have $\mathscr C^2=\deg \mathscr T=2m$, and we obtain (iii).

For (ii), let 
$$
r_L:H^0(Z,mF)^{G_1}\lras H^0(L, \mathscr O(2m)) 
$$
be the restriction map.
Then the restriction $\Phi_m^{G_1}|_L$ is exactly the rational map induced by the linear system corresponding to the image $r_L(H^0(mF)^{G_1})\subset H^0(\mathscr O(2m))$.
 Since $\Psi_m|_L:L\to\Lambda_m$ is a $2:1$ map, it has 2 branched points.
Let $z$ be an affine coordinate on $L$ such that $z=0$ and $z=\infty$ represent the branch points.
Then after a possible rescaling for the coordinate $z$, $\Psi_m|_L$ is written as $z\mapsto z^2$.
Hence as $\iota(\Lambda_m)$ is a rational normal curve in $\mathbb{CP}^m$, the composition $(\iota\circ\Psi_m)|_L$ is written as $z\mapsto (z^2,z^4,z^6,\cdots, z^{2m})$.
Then $\{1,z^2,z^4,\cdots,z^{2m}\}$ is a basis of the dual space of $r_L(H^0(mF)^G)$.

As is already proved, $\mathscr C=\Phi_m^{G_1}(L)$ is contained in a hyperplane in $\mathbb P(H^0(mF)^{G_1})^*$.
Moreover, $\mathscr C$ is non-degenerate in the hyperplane since $\Lambda_m$ is non-degenerate in $\mathbb P(H^0(Z,mF)^G)^*$ and $\pi_m|_{\mathscr C}$ is $2:1$. 
This means that 
\begin{align}\label{752}
\dim r_L(H^0(mF)^{G_1})=\dim r_L(H^0(mF)^G)+1.
\end{align}
Let $f(z)\in\mathbb C[z]$ be a polynomial such that  $\{1,z^2,z^4,z^6,\cdots,z^{2m},f(z)\}$ forms a basis of the dual space of $r_L(H^0(mF)^{G_1})$. 
As $r_L(H^0(mF)^{G_1})\subset H^0(\mathscr O_L(2m))$, we have $\deg f(z)\le 2m$.
Then by using a non-homogeneous coordinate on $\mathbb P^{\vee}H^0(mF)^{G_1}$, 
 $\Phi_m^{G_1}|_L$ can be written as 
\begin{align}\label{753}
z\longmapsto (z^2,z^4,z^6,\cdots,z^{2m}, f(z)).
\end{align}
Hence $\Phi_m^{G_1}|_L$ separates two points $z$ and $(-z)$ ($z\neq 0,\infty$) if and only if $f(z)\neq f(-z)$.
By \eqref{752}, $f(z)$ involves a term of odd-degree.
If $f(z)_{\rm{odd}}$ denotes the sum of all odd-degree terms, $f(z)= f(-z)$
holds iff $f(z)_{\rm{odd}}=0$.
If  $L$ is sufficiently general,  the equation 
$f(z)_{\rm{odd}}=0$ has 
exactly $2(m-1)$ solutions other than the obvious solution $z=0$.
This means that $\Phi_m^{G_1}|_L$ cannot  separate exactly $(m-1)$ pairs of points.
Namely $\Phi_m^{G_1}|_L$ identifies $(m-1)$ pairs of points. 
In order to show (ii) it remains to show that all the resulting singularities of $\mathscr C$ are nodes.
But this is obvious since the derivative of the map \eqref{753} has different values at $z$ and $(-z)$.
Thus we obtain (ii).

Finally, for (iv), since the map $\Phi_m^{G_1}$ preserves the real structure, $\mathscr C$ is real.
Further, as $L$ is non-singular and $\mathscr C$ is of course irreducible, the map $\Phi_m^{G_1}|_L:L\to \mathscr C$ is the normalization of the nodal curve $\mathscr C$.
Since $L$ has no real point, all the conditions $1^{\circ}$--$3^{\circ}$ in Section 4 are clearly satisfied.
\proofend

\begin{cor}\label{cor-existence}
For any $m\ge 2$, there exists a non-singular compact complex surface $S$ satisfying the following properties.
(i) $S$ has a real structure.
(ii) There exists a real rational curve $C$ on $S$  having $(m-1)$ nodes as its all singularities and satisfying $C^2=2m$.
(iii) The real structure satisfies the conditions $1^{\circ}, 2^{\circ}$ and $3^{\circ}$ in Section 4.
\end{cor}

\noindent Proof.
By taking a resolution of singularities of the complex surface $\mathscr T$ in Proposition \ref{prop-existence}, we only need to show that for any integer $m'\ge 2$ there exists a Joyce metric (or a $U(1)^2$-action) on $n\mathbb{CP}^2$ and a $U(1)^2$-invariant sphere such that 
the unique integer $m$ in Proposition \ref{prop:m} is exactly $m'$. 
But this is explicitly verified in the last part of \cite[\S 5.2]{Hon-mt}.
\proofend

\begin{rmk}{\em
It can be readily seen that in the proof of Corollary \ref{cor-existence} if we adapt the minimal resolution of the surface $\mathscr T$ then the surface $S$ is obtained from $\mathbb{CP}^1\times\mathbb{CP}^1$ by blowing-up $2m$ points satisfying the condition $(\ast)$ in the last subsection.
In particular the minimal resolution of $\mathscr T$ is a special form of the surface $S$ in Proposition \ref{prop-C*}.
Therefore, some of the surfaces in the last subsection admit a real structure and a real nodal rational curves.
However, the last real structure cannot be obtained as a natural lift of any real structure on $\mathbb{CP}^1\times\mathbb{CP}^1$, as the blowing-down map never preserves the real structure.
}\end{rmk}

\begin{rmk}{\em
As we have explained, the minitwistor spaces $\mathscr T$ in Proposition \ref{prop-existence} are
obtained as  quotient spaces of the twistor spaces of Joyce metrics on $n\mathbb{CP}^2$ by $\mathbb C^*$-action.
It might be worth mentioning that for  any one of these minitwistor spaces, there exist Moishezon twistor spaces on $n\mathbb{CP}^2$ with $\mathbb C^*$-action (for some $n$) whose quotient space is exactly the given minitwistor space but whose self-dual metric is {\em not} conformal to any Joyce metrics.
This is proved in \cite[Theorem 4.3]{Hon-mt}.
In the papers \cite{Hon07-3} and \cite{Hon07-4}, detailed structure  are studied for some of such twistor spaces.
}\end{rmk}

Finally we give a comment on the number of effective parameters involved in the construction of the surfaces $\mathscr T$ (which are the quotient spaces of the twistor spaces of the Joyce metrics).
As in the proof of Proposition \ref{prop-existence}, there is a canonical map $\pi_m|_{\mathscr T}:\mathscr T\to \iota(\Lambda_m)\simeq\mathbb{CP}^1$ whose fibers are conics.
For a general $\mathscr T$ (more precisely if $\mathscr T$ has no real singularities), the complex structure of $\mathscr T$ is uniquely determined by the discriminant locus of $\pi_m|_{\mathscr T}$, and the number of elements of the last locus is exactly $2m$.
Subtracting the dimension of the automorphism group of $\mathbb{CP}^1$,  the number of effective parameters is given by $2m-3$.
This coincides with the number obtained in the last subsection.
But of course, the location of discriminant locus is subject to a reality condition, and the last number is the dimension over real numbers.
On the other hand for the surface $\mathscr T$ which may have real singularities,  if $\nu_i$ denotes the number of real $A_i$-singularities (see \cite[Proposition 2.14]{Hon-mt}), then the number of reducible fibers of $\pi_m|_{\mathscr T}$ decreases, which is given by $2m-\sum_ii\nu_i$.



\small
\vspace{13mm}
\hspace{7.5cm}
$\begin{array}{l}
\mbox{Department of Mathematics}\\
\mbox{Graduate School of Science and Engineering}\\
\mbox{Tokyo Institute of Technology}\\
\mbox{2-12-1, O-okayama, Meguro, 152-8551, JAPAN}\\
\mbox{{\tt {honda@math.titech.ac.jp}}}\\
\mbox{{\tt {nakata@math.titech.ac.jp}}}\\
\end{array}$

\end{document}